\documentclass[nosumlimits]{amsart}
\usepackage{amsmath}
\usepackage{amssymb}

\newtheorem{thm}{Theorem}[section]
\newtheorem{defin}[thm]{Definition}
\newtheorem{lemma}[thm]{Lemma}
\newtheorem{cor}[thm]{Corollary}

\newtheorem{exa}{Example}
\newtheorem{pro}[thm]{Proposition}
\numberwithin{equation}{section}

\newcommand{\N}{\mathcal{N}}  
\newcommand{\Power}{\mathcal{P}}  

 \newcommand{\into}{\rightarrow}
 \newcommand{\impl}{\Longrightarrow}
 \newcommand{\less}{\setminus}
 \newcommand{\abs}[1]{\left\vert#1\right\vert}
 \newcommand{\set}[1]{\left\{#1\right\}}

 \newcommand{\norm}[1]{\left\Vert#1\right\Vert}
 
 \newcommand{\qtext}[1]{\quad\text{#1}\quad}
 \newcommand{\fa}{\qtext{for all}}
 \newcommand{\bb}{\begin{equation*}}
 \newcommand{\ee}{\end{equation*}}
 \newcommand{\bp}{\begin{proof}}
 \newcommand{\ep}{\end{proof}}

\begin{document}

\title[H.~Cartan's theorem and applications]
{On Henri Cartan's vectorial mean-value theorem\\ and its
applications\\ to Lipschitzian operators
and\\ generalized Lebesgue-Bochner-Stieltjes \\
 integration theory.
}
\author{ Victor M. Bogdan }

\address{Department of Mathematics, McMahon Hall 207, CUA, Washington DC 20064, USA}

\date{12 October 2009}

\email{bogdan@cua.edu}

\subjclass{46G10, 28A25, 47H99}
\keywords{Differentiability, vector analysis, vectorial calculus,
Dini derivatives, Lipschitzian functions, Lebesgue integral, Bochner integral}

\begin{abstract}
H. Cartan in his book on differential calculus proved a theorem
generalizing a Cauchy's mean-value theorem to the case of
functions taking values in a Banach space.

Cartan used this theorem in a masterful way to develop the entire
theory of differential calculus and  theory of
differential equations in finite and infinite dimensional Banach
spaces.

The author proves a generalization of this theorem to the case when
the inequality involving the derivatives holds everywhere with
exception  of a set of Lebesgue measure zero, and the
derivatives are replaced by  weaker derivatives. Namely the
right-sided Lipschitz derivative and lower right-sided Dini
derivative, respectively.

He also presents applications of the theorem to the study of
Lipschitzian operators in Banach spaces. Lipschitzian operators
played pivotal role in the n-body problems of electrodynamics, as
also in general n-body problem of Einstein's special theory of
relativity. For references see Bogdan\\
http://arxiv.org/abs/0909.5240 and\\
http://arxiv.org/abs/0910.0538.

Using the generalization of Cartan's theorem the author proves a version of
the fundamental theorem of calculus in a class of Bochner summable
functions. In the process he introduces the reader to the
generalized theory of Lebesgue-Bochner-Stieltjes integral and
Lebesgue and Bochner spaces of summable functions as developed by
Bogdanowicz.
\end{abstract}

\bigskip

\maketitle

\bigskip


\section{introduction}

\bigskip

Let $R$ and $Y$ denote, respectively, the space of reals
and a Banach space. We assume that the reader is
familiar with the notion of a Banach space as defined in \cite{banach},
\cite{cartan}, or \cite{DS1}.  Let $I=[a,b]$ be a closed bounded interval.
Henri Cartan in his book on differential calculus \cite{cartan},
proved the following

\bigskip

{\bf Theorem:} Let $I$ be the closed interval $[a,b].$
Assume that $f\!\!:\!I\mapsto Y$ and $g\!\!:\!I\mapsto R$ are
two continuous functions having right-sided derivatives $f'_r(x)$
and $g'_r(x)$ at every point of the open interval $(a,b).$ If
\begin{equation*}
 \|f'_r(x)\|\le g'_r(x)\fa x\in (a,b),
\end{equation*}
then
\begin{equation*}
 \|f(b)-f(a)\|\le g(b)-g(a).
\end{equation*}

\bigskip

Cartan used this theorem in
a masterful way to develop the entire theory of {\em Differential
Calculus} and {\em Theory of Differential Equations} in finite and infinite
dimensional Banach spaces.

We will prove a generalization of this
theorem to the case when the inequality involving
the derivatives holds everywhere with exception, perhaps,  of
a set of Lebesgue measure zero, and the derivatives are replaced by  weaker
derivatives. Namely the right-sided Lipschitz derivative
and lower right-sided Dini derivative, respectively.

We will also show applications of the theorem to the study of
Lipschitzian operators in Banach spaces. Lipschitzian operators
played pivotal role in the n-body problems of electrodynamics,
as also in general n-body problem of Einstein's special theory of relativity
\cite{einstein}.
For references see Bogdan \cite{bogdan61}--\cite{bogdan72}.

We will  show how one can use the theorem to prove
a version of the Fundamental Theorem of Calculus in a class
of Bochner summable functions.
In the process we will introduce the reader to the
generalized theory of Lebesgue-Bochner-Stieltjes
integral and Lebesgue and Bochner spaces of summable functions
as developed in Bogdanowicz \cite{bogdan10}--\cite{bogdan23}.

\section{Dini's one-sided derivatives}

\bigskip

Let $R$ denote the field of real numbers and $Y$ a Banach space.
We shall follow the notation used by Cartan \cite{cartan}.

\bigskip

\begin{defin}[Left- and right-sided derivatives]
A function $f:[a,b]\mapsto Y$ has a right-sided derivative at a
point $x\in [a,b)$ if \bb
    \lim_{h>0,h\into 0}\frac{1}{h}(f(x+h)-f(x))
\ee
exists. This limit will be denoted by $D_+f(x)=f'_r(x)$ and will be
called the right-sided derivative of $f$ at the point $x.$

Similarly we define
the left-sided derivative $D_-f(x)=f'_l(x)$ at a point $x\in (a,b]$
\bb
    D_-f(x)=\lim_{h<0,h\into 0}\frac{1}{h}(f(x+h)-f(x)).
\ee
Obviously a function $f$ has a derivative $Df(x)=f'(x)$ at
a point $x\in (a,b)$ if and
only if both derivatives $f'_r(x)$ and $f'_l(x)$ exist and
are equal. These derivatives represent an element from the
Banach space Y. In the case of reals, Y=R, it is sometimes
convenient to
admit also infinite values $\infty$ and $-\infty.$
\end{defin}

\bigskip

In some arguments it is convenient to introduce derivatives
known in the literature as Dini's derivatives.

\bigskip

\begin{defin}[Dini's derivatives]
Consider a real-valued function $f:(a,b)\mapsto R.$ By a
right-sided upper Dini derivative of $f$ at a point $x\in (a,b)$
we shall understand the  finite or infinite limit \bb
    D_+^u f(x)=\limsup_{h>0,h\into 0}\frac{1}{h}(f(x+h)-f(x)).
\ee
Similarly we define the right-sided lower Dini derivative by
\bb
    D_+^l f(x)=\liminf_{h>0,h\into 0}\frac{1}{h}(f(x+h)-f(x)).
\ee
And by analogy we define the left-sided upper derivative $D_-^u f(x)$
and the left-sided lower $D_-^l f(x)$ derivative. Clearly
The function $f$ has a right-sided derivative $D_+ f(x)$ at a point $x$
if and only if $D_+^u f(x)=D_+^l f(x).$
Similar relation is also valid for left-sided derivatives.
\end{defin}

\bigskip

\section{Lipschitzian functions and Lipschitzian derivatives}

\bigskip

\begin{defin}[Lipschitzian functions]Let $X,Y$ be some Banach spaces
and $I\subset X.$ We shall say that a function $f:I\mapsto Y$ is
Lipschitzian on the set $I$ if there exists a constant $m$ such
that \bb
    \|f(x_1)-f(x_2)\|\le m\|x_1-x_2\|\fa x_1,x_2\in I.
\ee

We shall say that such a function is Lipschitzian at a point $x_0\in I$
if there exist a constant $m<\infty$ and a positive $\delta$ such that
\bb
    \|f(x)-f(x_0)\|\le m \|x-x_0\|\fa x\in I, \|x-x_0\|<\delta.
\ee

When $X=R$ and $I=(a,b)$ is an open interval, we shall say that
the function is Lipschitzian at a  point $x_0$ to the right of it, if
for some $m<\infty$ there exists a positive $\delta$ such that
\bb
    \|f(x)-f(x_0)\|\le m|x-x_0|\fa x\in I, x_0<x<x_0+\delta.
\ee

Similarly we define what it means that the function is Lipschitzian
at a point $x_0$ to the the left of that point.
\end{defin}

\bigskip

\begin{defin}[Lipschitz derivatives]
Now let $Y$ denote a Banach space and consider a  function
$f:(a,b)\mapsto Y.$ By a right-sided Lipschitz derivative of $f$
at a point $x\in (a,b)$ we shall understand the  finite or
infinite limit \bb
    L_+ f(x)=\limsup_{h>0,h\into 0}\|\frac{1}{h}(f(x+h)-f(x))\|.
\ee
Similarly we define the left-sided Lipschitz derivative by
\bb
    L_- f(x)=\limsup_{h<0,h\into 0}\|\frac{1}{h}(f(x+h)-f(x))\|.
\ee
We shall say that the function $f$ has a Lipschitz derivative
$L f(x)$ at a point $x$
if  $L_+ f(x)=L_- f(x)$ and we denote the common value
by $L f(x).$
\end{defin}

\bigskip

Clearly if right-sided derivative $f'_r(x)$ exists then
the right-sided Lipschitz derivative exists and we have
the equality $\|f'_r(x)\|=L_+f (x).$ Similar relations
are valid for $f'_l(x)$ and $f'(x)$ and corresponding
Lipschitz derivatives $L_-f(x)$ and $Lf(x).$

Notice also that at a point $x\in (a,b)$ the Lipschitz
derivatives
\begin{equation*}
 Lf(x),\ L_+f(x),\ L_-f(x)
\end{equation*}
 are finite
if and only if the function is at the point $x$, respectively,
Lipschitzian,  Lipschitzian to the right,
Lipschitzian to the left of the point.

\bigskip

\section{Sets of Lebesgue measure zero}

\bigskip
\bigskip

\begin{defin}[Set of Lebesgue measure zero]
A set $A\subset R$ is said to be of {\bf Lebesgue measure zero} if
for every $\varepsilon>0$ there exists a countable collection
of intervals $T=\set{I_1,I_2,\ldots}$ such that
the set $A$ is contained in the union of sets in $T$
and
\begin{equation*}
 \sum_{I\in T} |I|\le \varepsilon,
\end{equation*}
where $\abs{I}$ denotes the length of the interval $I.$
\end{defin}

\bigskip

Clearly the empty set $\emptyset=(a,a),$ and any singleton $[b,b]$
is of Lebesgue measure zero. Moreover any countable set of points
forms a set of Lebesgue measure zero. Notice also that countable
union of sets of Lebesgue measure zero, is a set of Lebesgue
measure zero. Finally any subset of a set of Lebesgue measure zero
is a set of Lebesgue measure zero.

There exist also uncountable sets having Lebesgue measure zero.
A typical example of such a set  is the Cantor's set.
To construct Cantor's set take the closed interval $[0,1]$
and divide it into three equal intervals. From the middle remove
the open interval $(1/3,2/3).$ The remaining two closed intervals
have total length $2/3$ and they form a closed set $F_1$. Repeat
this process with each of the remaining intervals.

After
n-steps the remaining set $F_n$ will consist of the union of $2^n$
disjoint closed intervals of total length of $(2/3)^n.$
The sets $F_n$ are nested and their intersection $F$ will
represent a nonempty set of cardinality equal to cardinality of
the interval $[0,1].$ The set $F$ can be covered by a countable number
of intervals of total length as small as we please. Notice
that a finite cover by intervals we can always augment by
a sequence of intervals of the form $(a,a),$ that is by empty sets
to get a countable cover.

To prove that Cantor's set is of the same cardinality as the interval $[0,1]$
consider expansions into infinite fractions at the base $3$ of points belonging
to $F$
\begin{equation*}
 x=0.d_1,d_2,d_3,\ldots
\end{equation*}
where the digits $d_i\in\set{0,2}.$
Ignore the set of points which have  periodic expansions since
they represent some rational numbers that form a countable set.
Clearly all points that do not have periodic expansion are in the Cantor set $F.$

Similarly consider the binary expansions into nonperiodic
sequence of digits of points of the set $[0,1].$
\begin{equation*}
 y=0.a_1,a_2,a_3,\ldots
\end{equation*}
where $a_i\in\set{0,1}.$
Clearly the map $x\mapsto y$
given by the formula
\begin{equation*}
 a_i=d_i/2\fa i=1,2,\ldots
\end{equation*}
is one-to-one and onto. Thus the cardinalities of $F$
and $[0,1]$ are equal.

    It is clear that in the definition of a set $A\subset R$ of Lebesgue
    measure zero we can restrict ourself just to families consisting
    of open intervals. Indeed, take any $\varepsilon>0$ and let the sequence
    of intervals $I_n$ with end points $a_n,\,b_n$ be a cover of the set $A$
    with the total length of the intervals less than $\varepsilon/2.$
    Then the sequence
    \begin{equation*}
         J_n=(a_n-2^{-n-1}\varepsilon,b_n+2^{-n-1}\varepsilon)
    \end{equation*}
    consisting of open intervals will cover the set $T$ and its
    total length
    \begin{equation*}
     \sum_{n>0}J_n\le \sum_{n>0} I_n+
     \varepsilon\sum_{n>0} 2^{-n-1} \varepsilon\le\varepsilon.
    \end{equation*}
    Thus any set $A$ that can be covered by a sequence of interval of total
    length as small as we please can be covered by a sequence of open intervals
    of total length as small as we please. The converse is obvious.


 The following theorem characterizes the sets of Lebesgue measure
zero.

\bigskip

\begin{thm}[Characterization of sets of Lebesgue measure zero]
\label{null-set}
A set $A\subset R$ is of Lebesgue measure zero if and only if
there exists a sequence of open intervals $I_n$ such that
\begin{equation*}
 \sum_{n=1}^\infty \abs{I_n}\le 1
\end{equation*}
and
\begin{equation*}
 A\subset \bigcup_{n>k}I_n\fa k=1,2,\ldots
\end{equation*}
\end{thm}

\bigskip

\bp
%
    Assume that the set $A$ is of Lebesgue measure zero.
    Thus we can construct for every natural number $n$ a sequence of intervals
    \begin{equation*}
         J_{n,k}\,(k=1,2,\ldots)
    \end{equation*}
    such that
    \begin{equation*}
         A\subset \bigcup_{k=1}^\infty J_{n,k}\qtext{and}
         \sum_{k=1}^\infty |J_{n,k}|\le 2^{-n}.
    \end{equation*}
    Rearrange the double sequence $J_{n,k}$ into a single one $I_n$ and notice that
    it will have the properties stated in the theorem.

    \bigskip
    The converse of the above argument follows from the fact that
    remainder $r_n$ of a convergent series with terms $\abs{I_k}$
    converges to zero that is
    \begin{equation*}
         r_n=\sum_{k>n}\abs{I_k}\into 0\qtext{when}n\into\infty.
    \end{equation*}
\ep

\bigskip

\section{A theorem of Riesz}

\bigskip

The following theorem can be found
in the monograph of F. Riesz and B.~Sz.-Nagy \cite{riesz}.

\bigskip

\begin{thm}[Riesz]\label{riesz thm}
For every set $T\subset R$ of Lebesgue measure zero there exists
a continuous nondecreasing function $g\!:R\mapsto[0,1]$
such that
\begin{equation*}
 g'(x)=\infty\fa x\in T.
\end{equation*}
\end{thm}

\bigskip

\bp
    For any open interval $I=(a,b)$ define the function $G_I:R\mapsto R$
    by the formula
    \begin{equation}\label{function G sub I}
     \begin{array}{llc}
    G_I(x)=0   &\text{ if }&x<a,\\
    G_I(x)=x-a &\text{ if }&a\le x\le b,\\
    G_I(x)=b-a &\text{ if }&x>b.
    \end{array}
    \end{equation}

    For the set $T$ let $I_n$ be a sequence of open
    intervals as in Theorem \ref{null-set}.
    Let $g_n=G_{I_n}$ for all $n.$ Consider the series with terms $g_n.$
    It consists of nonnegative nondecreasing continuous functions. Since
    \begin{equation*}
         g_n(x)\le \abs{I_n}\fa x\in R,\ n=1,2,\ldots
    \end{equation*}
    the series converges uniformly on the entire space $R$ to a continuous function.
    Thus the function $g$ given by the formula
    \begin{equation*}
     g(x)=\sum_{n=1}^\infty g_n(x)
    \end{equation*}
    is well defined and represents a nonnegative nondecreasing and continuous function.

    Now let us prove that at every point of the set $T$ the derivative of the function
    $g$ is equal to infinity.
    To this end notice that every point $t\in T$ belongs to an infinite
    number of the intervals $I_n.$ Let $k_n$ denote the number of intervals in the
    sequence of the first $n$ intervals with index $\le n,$
    \begin{equation*}
         I_1, I_2, \ldots, I_n
    \end{equation*}
    containing the point $t.$ Let $J_n$ denote the intersection of these intervals.
    Since $t$ belongs to all of these intervals and the intersection of a finite
    number of open sets is open, the set $J_n$ is open. So there is a positive
    number $\delta>0$ such that
    \begin{equation*}
     (t-\delta,t+\delta)\subset J_n.
    \end{equation*}
    Now consider the difference quotient
    \begin{equation*}
     \frac{g(x)-g(t)}{x-t}
    \end{equation*}
    for $|x-t|<\delta.$ We have the following representation
    \begin{equation*}
     \frac{g(x)-g(t)}{x-t}=\sum_{j\le n}\frac{g_j(x)-g_j(t)}{x-t}+
    \sum_{j> n}\frac{g_j(x)-g_j(t)}{x-t}
    \end{equation*}
    Since each function $g_j$ is nondecreasing each term of each sum is nonnegative.
    In the first sum there will be at least $k_n$ terms for which the
    corresponding difference quotient will be equal to $1.$ Since the second
    sum consists of nonnegative terms we have the lower estimate
    \begin{equation*}
     \frac{g(x)-g(t)}{x-t}\ge k_n\fa |x-t|<\delta.
    \end{equation*}
    Since $k_n\into\infty$ the above estimate proves that $g'(t)=\infty.$
\ep

\bigskip

In the sequel we will need the following lemma.

\bigskip

\begin{lemma}\label{psi function}
 Let $I=[a,b]$ denote a closed bounded interval.
Let $T$ be a  set of Lebesgue measure zero lying in its interior
$(a,b).$ There exists a positive, increasing, continuous function
$\psi:I\mapsto R$ such that
\begin{equation*}
    \psi'(t)=\infty \quad\text{at every point}\quad t\in T.
\end{equation*}
\end{lemma}

\bigskip

\bp
    Take the function $g$ as in the Riesz theorem
    (\ref{riesz thm}) corresponding to the
    the set $T,$ and the function $G_I$ as defined in (\ref{function G sub I}).
    Define the function $\psi$ by the formula
\begin{equation*}
     \psi(x)=g(x)+G_I(x)+1\fa x\in [a,b].
\end{equation*}
    Clearly this function will satisfy the conditions of the lemma.
\ep

\bigskip

\section{Strong Mean-Value theorem}

\bigskip

We shall prove the following theorem representing a generalization
of a theorem due to Henri Cartan \cite{cartan},  p. 37.

\begin{thm}[Strong Mean-Value Theorem]
Let $I=[a,b]$ be a closed bounded interval and let $T\subset(a,b)$ be a set of
Lebesgue measure zero. Let $f:I\mapsto Y$ and $g:I\mapsto R$ be
continuous functions. Assume that the right-sided Lipschitz
derivative $L_+f(x)$ and right-sided lower Dini derivative $D_+^l
g(x)$ exist and are finite at every point $x\in (a,b).$ Then the
inequality
\begin{equation}\label{cond on derivatives}
    L_+ f(x)\le D_+^l g(x)\fa x\in (a,b)\less T,
\end{equation}
implies the inequality
\begin{equation}\label{str mean value}
    \|f(b)-f(a)\|\le g(b)-g(a).
\end{equation}
\end{thm}

\bigskip

\bp
    Let $\psi$ denote the function from the lemma \ref{psi function}
    corresponding
    to the set $T.$

    To prove the
    theorem it is sufficient to prove the following.
    For every $\varepsilon>0$ we have
    \begin{equation}\label{cond A}
        \|f(x)-f(a)\|\le
        g(x)-g(a)+\varepsilon((x-a)+\psi(x))\fa x\in I.
    \end{equation}
    Once the validity of the statement is established
    setting $x=b$ and passing
    to the limit $\varepsilon\into 0$ will yield the inequality
    (\ref{str mean value}).

    We will prove the validity of the above
    condition by contradiction.
    Assuming that the statement is not true
    we get that for some $\varepsilon>0$ the set
    \begin{equation}\label{not A}
        U=\set{x\in I:\ \|f(x)-f(a)\|>
        g(x)-g(a)+\varepsilon((x-a)+\psi(x))}
    \end{equation}
    is nonempty. Thus from the axiom of continuity follows that
    the number $c=\inf U$ is well defined. Since all the functions
    in the inequality (\ref{cond A}) are continuous, taking all
    functions from the righthand side  of
    the inequality onto the left side and
    denoting the function on the left side by $\phi,$ we get a
    representation
    of the set $U$ in the form
    \bb
        U=\set{x\in I:\,\phi(x)>0}=\phi^{-1}(0,\infty),
    \ee
    where $\phi$ is a continuous function. Thus $U$ is an open set
    in $I$ as the inverse image of an open set by means of a continuous
    function.

    Notice that at $x=a$ the inequality (\ref{cond A})
    is strict since $\psi(0)>0.$
    Thus it follows from continuity of the functions involved that for
    some $d>a$ the inequality (\ref{cond A}) holds for all $x\in [a,d).$
    This means that $U\subset[d,b]$ and thus $a<d\le c,$ and so $c\not=a.$

     We also have  $c\not=b.$ Otherwise $U=\set{b}$
     and the set $U$ would not be
     open in $I.$ Hence we must have that $a<c<b.$

    The set $U$ cannot contain the point $c.$
    Otherwise some closed interval
    $[c-\delta,c+\delta]\subset (a,b),$ where $\delta>0,$
    would be in the set $U.$
    Thus $c>c-\delta,$ that is the greatest lower bound $c$ of  the
    set $U$ is grater than an element $c-\delta$ of the set $U$
    yielding $c>c-\delta\ge c,$ a contradiction.

    Since $c\not\in U$ we must have
    \begin{equation}\label{B0}
        \|f(c)-f(a)\|\le g(c)-g(a)+\varepsilon((c-a)+\psi(c)).
    \end{equation}

    Put $\eta=L_+ f(c)$ and $\tau=D_+^l g(c).$ From the definitions
    of $\limsup$ and $\liminf$ we get
    \begin{equation}\label{A1}
        \eta=\inf_{\delta>0}\ \sup
        \set{0<h<\delta:\ \frac{1}{h}\|f(c+h)-f(c)\|}
    \end{equation}
    and
    \begin{equation}\label{A2}
        \tau=\sup_{\delta>0}\ \inf
        \set{0<h<\delta:\ \frac{1}{h}[g(c+h)-g(c)]}.
    \end{equation}

    From (\ref{A1}) follows that for $\eta+\varepsilon/2$ there
    exists a $\delta_1>0$ such that
    \begin{equation}\label{A3}
        \frac{1}{h}\|f(c+h)-f(c)\|<\eta+\varepsilon/2
        \qtext{if} 0<h<\delta_1.
    \end{equation}
    Similarly from (\ref{A2}) we get that there exists $\delta_2>0$
    such that
    \begin{equation}\label{A4}
        \tau-\varepsilon/2< \frac{1}{h}[g(c+h)-g(c)]
        \qtext{if} 0<h<\delta_2.
    \end{equation}

    \bigskip
    We have two possibilities that remain either $c\not\in T$ or $c\in T.$
    \bigskip

    First let us consider the case when $c\not\in T.$
    In this case by assumption of the theorem $\eta\le\tau.$
    The conditions (\ref{A3}) and (\ref{A4}) yield that for
    some positive $\delta<\min\set{\delta_1,\delta_2}$ we have
    \begin{equation}\label{A5}
        \frac{1}{h}\|f(c+h)-f(c)\|\le\eta+\varepsilon/2
        \le \tau+\varepsilon/2\le\frac{1}{h}[g(c+h)-g(c)]+\varepsilon
        \qtext{if} 0<h\le\delta.
    \end{equation}
    Equivalently
    \begin{equation}\label{A6}
        \|f(x)-f(c)\|
        \le [g(x)-g(c)]+\varepsilon(x-c)
        \qtext{if} c<x\le c+\delta.
    \end{equation}
    Thus from the relations (\ref{B0}) and (\ref{A6}), and the fact
    that the function $\psi$ is increasing, follows that
    \begin{equation}\label{A7}
        \begin{split}
        \|f(x)-f(a)\|&\le \|f(x)-f(c)\|+\|f(c)-f(a)\| \\
        &\le
        [g(x)-g(a)]+\varepsilon[(x-a)+\psi(x)]
        \qtext{if} c<x\le c+\delta.\\
        \end{split}
    \end{equation}
    Thus the number $c+\delta$ is a lower bound of the set $U$
    greater then the greatest lower bound $c=\inf U.$ A contradiction.

    Finally let us consider the case $c\in T.$ Since $\psi'(c)=\infty,$
    we get that for the number $\eta-\tau+\varepsilon$ there exists
    a positive number $\delta_3$ such that
    \begin{equation}\label{B1}
        \eta-\tau+\varepsilon<\varepsilon \frac{1}{h}(\psi(c+h)-\psi(c))
        \qtext{if} 0<|h|<\delta_3.
    \end{equation}\bb\ee
    Thus we get
    \begin{equation}\label{B2}
        \eta+\varepsilon/2<(\tau-\varepsilon/2)
        +\varepsilon \frac{1}{h}(\psi(c+h)-\psi(c))
        \qtext{if} 0<|h|<\delta_3.
    \end{equation}

    Hence for a positive $\delta<\min\set{\delta_1,\delta_2,\delta_3 }$
    from the relations (\ref{A3}) and (\ref{A4}) we get
    \begin{equation}\label{B3}
        \frac{1}{h}\|f(c+h)-f(c)\|\le\frac{1}{h}[g(c+h)-g(c)]+
        \varepsilon\frac{1}{h}[\psi(c+h)-\psi(c)]
        \qtext{if} 0<h\le \delta
    \end{equation}
    and therefore
    \begin{equation}\label{B4}
        \|f(x)-f(c)\|\le[g(x)-g(c)]+
        \varepsilon[\psi(x)-\psi(c)]
        \qtext{if} c<x\le c+\delta.
    \end{equation}
    Thus from the relations (\ref{B0}) and (\ref{B4}) we get the relation
    \bb
        \begin{split}
        \|f(x)-f(a)\|&\le \|f(x)-f(c)\|+\|f(c)-f(a)\| \\
        &\le
        [g(x)-g(a)]+\varepsilon[(x-a)+\psi(x)]
        \qtext{if} c<x\le c+\delta.\\
        \end{split}
    \ee
    and this relation as in the previous case leads
    to a contradiction.
\ep

\bigskip

\section{Henri Cartan's Mean-Value Theorem}

\bigskip

\begin{cor}[Cartan]
Let $I=[a,b]$ be a closed bounded interval. Let $f:I\mapsto Y$ and
$g:I\mapsto R$ be continuous functions. Assume that the
right-sided derivatives $f'_r(x)$ and $g'_r(x)$ exist at every
point $x\in (a,b).$ Then the inequality \bb
     \|f'_r(x)\|\le g'_r(x)\fa x\in (a,b),
\ee
implies the inequality
\bb\label{Cartan mean est}
    \|f(b)-f(a)\|\le g(b)-g(a).
\ee
\end{cor}

\bigskip

\bp
    To prove the theorem notice that the existence of
    right-sided derivatives $f'_r(x)$ and $g'_r(x)$
    implies the existence and finiteness of the derivatives
    $L_+f(x)=\|f'_r(x)\|$ and $D_+^lg(x)=g'_r(x).$ Thus
    we can use the previous theorem.
\ep

\bigskip

\begin{cor}
It is obvious that there exist corresponding theorems for
left-sided derivatives. Is is sufficient to introduce functions
$\tilde{f}(x)=f(-x)$ and $\tilde{g}(x)=g(-x)$ and observe
that right-sided derivatives are mapped onto left-sided derivatives.
\end{cor}

\section{A condition for a function to be nondecreasing}

\bigskip

\begin{cor}
Let $g:[a,b]\mapsto R$ be a continuous function whose right-sided
lower Dini derivative $D_+^lg(x)$ is finite for all $x$ in
$(a,b).$ If the derivative is nonnegative at every point of the
interval $(a,b),$ with exception perhaps of a set of Lebesgue
measure zero, then the function $g$ is nondecreasing on the entire
interval $[a,b].$
\end{cor}

\bigskip

\bp
    Define $f(x)=0$ for all $x\in [a,b].$  Clearly $L_+f(x)=0$ on $(a,b).$
    Assume $x_1<x_2$ are any two fixed points of the interval $[a,b].$
    If you consider the function $g$ on the interval $[x_1,x_2]$ then
    all conditions of the Strong Mean-Value Theorem are satisfied on
    the interval $[x_1,x_2].$ Hence
    \bb
        0=\|f(x_2)-f(x_1)\|\le g(x_2)-g(x_1)\fa x_1<x_2.
    \ee
    This shows that the function $g$ is nondecreasing.
\ep

\bigskip

\section{A sufficient condition for a vector function to be constant}
\bigskip

\begin{cor}[Sufficient condition for constancy of a vector function]
\label{cor on constancy}
Let $J$ be an open interval and $Y$ a Banach space.
 Assume that $f:J\mapsto Y$ is continuous and $f'_r(t)=0$
 for  almost all $t\in J.$

Then the  function $f$ is constant on the interval $J.$
\end{cor}

\bigskip
\bp
    Introduce function $g(t)=0$ for all $t\in J.$
    Take any two different points $t_1,t_2\in J.$ We may assume
    that $t_1<t_2.$ Let $f(t_1)=y_0.$
    The pair of functions $f,g$ considered on the interval
    $[t_1,t_2]$ satisfies the assumptions of the Strong Mean-Value
    Theorem. Thus we have
    \begin{equation*}
     \norm{f(t_2)-f(t_1)}\le g(t_2)-g(t_1)=0
    \end{equation*}
    That is $f(t)=y_0$ for all $t\in J.$
\ep

\bigskip

\section{Characterization of Lipschitzian vectorial functions on intervals}

\bigskip

\begin{cor} Assume that $Y$ is a Banach space and $J$ any interval.
Let $f:J\mapsto Y$ be a continuous function whose right-sided
Lipschitz derivative $L_+f(x)$ is finite for all $x$ in $(a,b).$
If there exists a constant $m<\infty$ such that \bb
    L_+f(x)\le m
\ee
at every point of the interval $(a,b),$ with exception perhaps
of a set of Lebesgue measure zero, then the
function $f$ is Lipschitzian on the entire interval $[a,b]$ and
\bb
    \|f(x)-f(y)\|\le m |x-y|\fa x,y\in [a,b].
\ee
The converse condition is obvious.
\end{cor}

\bigskip

\bp
    Take any two numbers $x,y$ in the interval $[a,b].$ Without
    loss of generality we may assume $x>y.$
    Put $g(x)=mx$ for all $x\in [a,b].$ Consider the closed
    interval $[y,x].$ It follows from  the Strong Mean-Value Theorem
    that
    \begin{equation*}
        \|f(x)-f(y)\|\le m x-m y=m|x-y|\fa x,y\in [a,b], x>y.
    \end{equation*}
\ep

\bigskip

\section{Characterization of Lipschitzian functions on convex sets}

\bigskip

\begin{defin}[Lipschitz derivative]
Now let $X,Y$ denote Banach spaces and $W$ a convex set in $X.$
Consider a  function $f:W\mapsto Y.$
By a  {\bf Lipschitz derivative} of $f$ at a point
$x\in W$ we shall understand the  finite
or infinite limit
\bb
    L f(x)=\limsup_{y\into x,y\neq x}\frac{\norm{f(y)-f(x)}}{\norm{y-x}}.
\ee
\end{defin}

\bigskip

\begin{thm} Assume that $X$ and $Y$ are  Banach spaces
and $W$ is a convex set in $X.$

Let $f:W\mapsto Y$ be a continuous function whose
Lipschitz derivative $Lf(x)$ is finite
for all $x$ in $W.$ If there exists a constant $m<\infty$
such that
\begin{equation*}
    Lf(x)\le m
\end{equation*}
at every point of the set $W,$ with exception perhaps
of a countable set, then the
function $f$ is Lipschitzian on the set $W$ and
\bb
    \norm{f(x)-f(y)}\le m \norm{x-y}\fa x,y\in W.
\ee
The converse is obvious.
\end{thm}

\bigskip
\section{Right-sided antiderivative of a function}

\bigskip
\begin{defin}[Right-sided antiderivative]
Assume that $J$ is an open interval in $R$ and $Y$ is a Banach
space. Assume that $f,g$ are two functions defined on $J$
into the Banach space $Y.$
If the function $f$ is continuous and has
right-sided derivative $f'_r(t)$ at every point $t\in J$
and for some set $T$ of Lebesgue measure zero we have the equality
\begin{equation*}
 f'_r(t)=g(t)\fa t\in J\less T,
\end{equation*}
then we shall say the the function $f$ forms a {\bf right-sided antiderivative} of
the function $g$ over the interval $J.$
\end{defin}

\bigskip
Notice that any continuous piecewise linear function or more generally any
continuous piecewise differentiable function forms
a right-sided antiderivative of its
derivative in the above sense. So it is an essential generalization
of the notion of antiderivative.

This generalization has another important
property stated in the following theorem.

\bigskip
\begin{thm}[Any two right-sided antiderivatives
of the same function differ by a constat]\label{two anti-d differ byconstant}
Assume that $J$ is an open interval and $Y$ a Banach space.
Assume that $f_1,f_2,g$ are defined on the interval
$J$ into the Banach space $Y.$
If both $f_1$ and $f_2$ form right-sided derivatives of the function $g$ over
the interval $J$
then there exists a vector $y_0\in Y$ such that
\begin{equation*}
 f_1(t)-f_2(t)=y_0\fa t\in J.
\end{equation*}
\end{thm}

\bigskip
\bp
    Put $f=f_1-f_2.$ The function $f$ is continuous and from the linearity
    of the limit operation we get that
    \begin{equation*}
    f'_r(t)=D_+f(t)=D_+f_1(t) - D_+f_2(t)\fa t\in J
    \end{equation*}
    Moreover
    \begin{equation*}
    f'_r(t)=0\qtext{for almost all}t\in J.
    \end{equation*}
    Using the Corollary (\ref{cor on constancy}) we can conclude
    that for some vector $y_0\in Y$ we have
    \begin{equation*}
    f_1(t)-f_2(t)=y_0\fa t\in J.
    \end{equation*}
\ep

\bigskip

Now a natural question arises when we can use the
formula known as the Fundamental Theorem of Calculus
\begin{equation*}\label{fun thm calculus}
 \int_{t_1}^{t_2}f'_r(u)\,du=f(t_2)-f(t_1)\fa t_1,t_2\in J.
\end{equation*}

We shall prove later that if the function $f'_r$ is locally bounded on the set
$J,$ excluding perhaps a set of Lebesgue measure zero, then the function $f'_r$
is Bochner summable on every closed interval $[t_1,t_2]$
and the the formula (\ref{fun thm calculus}) holds.

\bigskip


\section{Vectorial Lebesgue-Bochner Integration Theory}

\bigskip

In this section we shall present a development of the theory of
Lebesgue and Bochner spaces of summable functions and present a
construction and fundamental theorems of the theory. We shall
follow the approach of Bogdanowicz \cite{bogdan10} and
\cite{bogdan14} with some modifications. In the process we shall
construct a generalized Lebesgue-Bochner-Stieltjes integral as
developed \cite{bogdan10}.

The development of the integration theory beyond the classical
Riemann integral is essential in the modern theory of differential
equations, theory of operators, probability, and optimal control
theory, and most important in theoretical physics.

\bigskip

Assume that $Y, Z, W$ represent some Banach spaces either over the
field $R$ of reals or over the field $C$ of complex numbers.

Denote by $U$ the space of all bilinear bounded operators $u$ from
the space $Y\times Z$ into $W$. Norms of elements in the spaces
$Y, Z, W, U$ will be denoted by $|\ \,|$.

If $V$ is any nonempty collection of subsets of an abstract space
$X$ denote by $S(V)$ the family of all sets that are disjoint
unions of finite collections of sets from the collection $V.$
Since the empty collection is finite, we implicitly assume that
the empty set $\emptyset$ belongs to $S(V).$ The family $S(V)$
will be called the {\bf family of simple sets} generated by the
family $V.$

A nonempty family of sets $V$ of the space $X$ is called a {\bf
prering} if the following conditions are satisfied: if $A_{1},
A_{2}\in V$, then both the intersection $ A_{1}\cap A_{2}$ and set
difference $ A_{1}\less A_{2}$
belong to the family $S(V).$ 

A family of sets $V$ of the space $X$ is called a {\bf ring} if
$V$ is a prering such that $V=S(V)$ which is equivalent to the
following conditions: if $A_{1}, A_{2}\in V$, then $ A_{1}\cup
A_{2}\in V$ and $ A_{1}\less A_{2}\in V$. It is easy to prove that
a family $V$ forms a prering if and only if the family $S(V)$ of
the simple sets forms a ring. Every ring (prering) $V$ of a space
$X$ containing the space $X$ itself is called an {\bf algebra
(pre-algebra)}, respectively.

If the ring $V$ is  closed under countable unions it is called a
{\bf sigma ring} ($\sigma$-ring for short.) If the ring $V$ is
closed under countable intersections it is called a {\bf delta
ring} ($\delta$-ring for short.) It follows from de Morgan law
that $\delta$-algebra and $\sigma$-algebra represent the same
notion.

A finite-valued function $v$ from a prering $V$ into $\langle
0,\infty),$ the non-negative reals, satisfying the following
implication
\begin{equation}\label{countable additivity}
    A=\bigcup_{t\in T}A_t\impl v(A)=\sum_{t\in T}v(A_t)
\end{equation}
for every set $A\in V,$ that can be decomposed into countable
collection $A_t\in V\,(t\in T)$ of disjoint sets, will be called a
{\bf $\sigma$-additive positive measure.} It was called a {\bf
positive volume} in the preceding papers of the author. Notice
that since by definition every prering $V$ contains at least one
element $A\in V$, we must have that $\emptyset=A\less A\in V.$
Thus from countable additivity (\ref{countable additivity})
follows that $v(\emptyset)=0.$

By {\bf Lebesgue measure} over an abstract space $X$ we shall
understand any set function $v$ from a $\sigma$-ring $V$ of the
space $X$ into the extended non-negative reals $\langle
0,\infty\rangle,$ that satisfies the implication (\ref{countable
additivity}) and has value zero on the empty set $v(\emptyset)=0.$
We have to postulate this explicitly to avoid the case of a
measure that is identically equal to $\infty.$

As in Halmos \cite{halmos} a triple $(X,V,v),$ where $X$ denotes
an abstract space and $V$  a prering of the space $X$ and $v$ a
$\sigma$-additive non-negative finite-valued measure on the
prering $V,$ will be called a {\bf measure space.}

Halmos considered such measure spaces for the case when $V$ forms
a ring of sets. Since every ring satisfies the axioms of a
prering, our notion of a measure space is more general.

Halmos used such measure spaces to construct Lebesgue measures and
to base on them the development of  the integration theory. We
reverse the process by first developing the integration theory and
obtain the Lebesgue measure as a by product.

\bigskip

It is clear that every finite Lebesgue measure forms a positive
measure in our sense, and in the case when it has infinite values
by striping it of infinities we obtain a positive measure.

The development of the classical Lebesgue-Bochner theory of the
integral goes through the following main stages as in Halmos
\cite{halmos} and Dunford and Schwartz \cite{DS1}:
\begin{itemize}
\item The construction and development of the Caratheodory theory
of outer measure $v^*$
        over an abstract space $X.$
\item The construction of the Lebesgue measure $v$ on the sigma
ring $V$ of measurable sets
        of the space $X$ induced by the outer measure $v^*.$
\item The development of the theory of real-valued measurable
functions $M(v,R).$ \item The construction of the Lebesgue
integral $\int f\,dv.$ \item The construction and development of
the theory of the space $L(v,R)$
        of Lebesgue summable functions.
\item The construction and development of the theory of the space
$M(v,Y)$ of
        Bochner measurable functions.
\item The construction of the Bochner integral $\int f\,dv$ and of
the space
        $L(v,Y)$ of Bochner summable functions $f$ from the space $X$ into
        any Banach space $Y.$
\end{itemize}

The construction of the Lebesgue's integral is an abstraction from
the area under the graph of the function similar to the ideas of
Riemann though different in execution.

From the point of view of Functional Analysis both the Lebesgue
and Bochner integrals are particular linear continuous operators
from the space $L(v,Y)$ of summable functions into the Banach
space $Y.$ Moreover from the theory of the space $L(v,Y)$ one can
easily derive the theory of the spaces $M(v,Y),$ $M(v,R),$ and
$L(v,R)$ and of the Lebesgue and Bochner integrals and also the
theory of Lebesgue measure. For details see Bogdanowicz
\cite{bogdan14} and \cite{bogdan23}.

We shall show in brief how one can develop the theory of the space
$L(v,Y)$ and to construct an integral of the form $\int
u(f,d\mu),$ where $u$ is any bilinear operator from the product
$Y\times Z$ of Banach spaces into a Banach space $W$ and $\mu$
represents a vector measure. This integral for the case, when the
spaces $Y,Z,W$ are equal to the space $R$ of reals and the
bilinear operator $u$ represents multiplication $u(y,z)=yz$,
coincides with the Lebesgue integral
\begin{equation}\label{Lebesgue integral}
    \int f\,dv=\int u(f,dv)\fa f\in L(v,R).
\end{equation}
In the case, when $Y=W$ and $Z=R$ and $u(y,z)=zy$ represents the
scalar multiplication, the integral coincides with the Bochner
integral
\begin{equation}\label{Bochner integral}
    \int f\,dv=\int u(f,dv)\fa f\in L(v,Y).
\end{equation}

\bigskip

It is good to have a few examples of the measure spaces. The first
example corresponds to Dirac's $\delta$ function.

\bigskip

\begin{exa}[{\bf Dirac measure space}] Let $X$ be any abstract set and $V$
the family of all subsets of the space $X.$ Let $x_0$ be a fixed
point of $X.$ Let $v_{x_0}(A)=1$ if $x_0\in A$ and  $v_{x_0}(A)=0$
otherwise. Since $V$ forms a sigma ring the triple $(X,V,v)$ forms
in this case a Lebesgue measure space.
\end{exa}

\bigskip

\begin{exa}[{\bf Counting measure space}] Let $X$ be any abstract set and
$V=\set{\emptyset,\set{x}:\ x\in X}.$ Let $v(A)=1$ for all
singleton sets $A=\set{x}$ and $v(\emptyset)=0.$ The triple
$(X,V,v)$ forms a measure space that is not a Lebesgue measure
space.
\end{exa}

\bigskip

\begin{exa}[{\bf Striped Lebesgue measure space}]
Assume that $M$ is a sigma-ring of subsets of $X$ and $\mu$  is
any Lebesgue measure on $M.$ Let $$V=\set{A\in
M:\,\mu(A)<\infty}.$$ Plainly $V$ forms a ring and thus a prering.
Then the restriction $\mu$ to $V$ yields a measure space
$(X,V,\mu.)$
\end{exa}

\bigskip

The most important measure space to the sequel is the following.

\begin{pro}[{\bf Riemann measure space}]
\label{Riemann measure space} Let $R$ denote the space of reals
and $V$ the collection of all bounded intervals $I$ open, closed,
or half-open. If $a\le b$ are the end points of an interval $I$
let $v(I)=b-a.$ Then the triple $(R,V,v)$ forms a measure space.
We shall call this space the {\bf Riemann measure space.}
\end{pro}

\bigskip

\bp
    The collection $V$ of intervals forms a prering. Indeed the intersection
    of any two intervals is an interval or an empty set. But empty set
    can be represented as an open interval $(a,a)=\emptyset.$
    The set difference of two intervals is either the union of two disjoint intervals
    or a single interval or an empty set. Thus we have that for any two intervals
    $I_1,\,I_2\in V$ we have $I_1\cap I_2\in S(V)$ and $I_1\less I_2\in S(V).$
    This proves that $V$ is a prering.

    To prove countable additivity
    assume that we have a decomposition of an interval $I$ with ends $a\le b$
    into disjoint countable collection $I_t(t\in T)$ of intervals with end
    points $a_t\le b_t,$ that is
    \begin{equation}\label{decomposition}
        I=\bigcup_{t\in T}I_t.
    \end{equation}
     The case when interval $I$ is empty or consists
    of a single point is obvious. So without loss of generality we may assume
    that the interval $I$ has a positive length and that
    our index set $T=\set{1,2,3,\ldots}.$ Take any $\varepsilon>0$ such that
    $2\varepsilon<v(I).$ Let $I^\varepsilon=[a+\varepsilon,b-\varepsilon]$ and
    $I^\varepsilon_t=(a_t-\varepsilon2^{-t},b_t+\varepsilon2^{-t})$ for all $t\in T.$

    The family $I^\varepsilon_t(t\in T)$ forms an open cover of the compact
    interval $ I^\varepsilon$ thus there exists a finite set $J\subset T$ of indexes
    such that
    \begin{equation*}
        I^\varepsilon\subset \bigcup_{t\in J} I^\varepsilon_t.
    \end{equation*}
    The above implies
    \begin{equation*}
        \begin{split}
        v(I)-2\varepsilon&=v(I^\varepsilon)\le \sum_{t\in J}v(I^\varepsilon_t)
        \le \sum_{t\in T}v(I^\varepsilon_t)\\
        &=\sum_{t\in T}(v(I_t)+2^{-t+1}\varepsilon)=
        \sum_{t\in T}v(I_t)+2\varepsilon.
        \end{split}
    \end{equation*}
    Passing to the limit in the above inequality when $\varepsilon\into 0$
    we get
    \begin{equation*}
        v(I)\le \sum_{t\in T}v(I_t)
    \end{equation*}
    On the other hand from the relation (\ref{decomposition})  follows that
    for any finite set $J$ of indexes we have
    \begin{equation*}
        I\supset \bigcup_{t\in J}I_t\impl v(I)\ge \sum_{t\in J}v(I_t).
    \end{equation*}
    Since $\sup_J \sum_{t\in J}v(I_t)=\sum_{t\in T}v(I_t)$ we get from
    the above relations that the set function $v$ is countably additive
    and thus it forms a measure.
\ep

\bigskip

As will follow from the development of this theory the Riemann
measure space generates the same space of summable functions and
the integral as the classical Lebesgue measure over the reals. It
is good to see a few more examples of measures related to this
one.

\bigskip

\begin{pro}[{\bf Stieltjes measure space}]
Let $R$ denote the space of reals, and $g$  a nondecreasing
function from $R$ into $R,$ and $D$ the set of its discontinuity
points of $g.$ Let $V$ denote the collection of all bounded
intervals $I$ open, closed, or half-open with end points
$a,b\not\in D.$ If $a\le b$ are the end points of an interval $I$
let $v(I)=f(b)-f(a).$ Then the triple $(R,V,v)$ forms a measure
space.
\end{pro}

\bigskip

\bp
    The proof is similar to the preceding one and we leave it to the reader.
\ep

\bigskip

A nondecreasing left-side continuous function $F$ from the
extended closed interval $E=\langle-\infty,+\infty\rangle$ such
that $F(-\infty)=0$ and $F(+\infty)=1$ is called a {\bf
probability distribution function.} Any measure space $(X,V,v)$
over a prering $V$ such that $X\in V$ and $v(X)=1$ is called a
{\bf probability measure space.}

\bigskip

\begin{pro}[Probability distribution generates probability measure space]
Let $F$ be a probability distribution on the extended reals $E.$
Let $V$ consists of all intervals of the form $\langle a,b)$ or
$\langle a,\infty\rangle,$ where $a,b\in E.$ If $I\in V$ let
$v(I)$ denote the increment of the function on the interval $I$
similarly a in the case of Stieltjes measure space.

Then the triple $(E,V,v)$ forms a probability measure space.
\end{pro}

\bigskip

\bp
    To prove this proposition notice that the space $E$ can be considered
    as compact space and the proof can proceed similarly as in the case of the
    Riemann measure space.
\ep

\bigskip

In the case of topological spaces there are two natural prerings
available to construct a measure space: The prering consisting of
differences $G_1\less G_2$ of open sets, and the prering
consisting of differences $Q_1\less Q_2$ of compact sets.

\bigskip

\section{Construction of the elementary integral spaces}

\bigskip

\begin{defin}[Vector measure]
A set function $\mu$ from a prering $V$ into a Banach space $Z$ is
called a {\bf vector measure} if for every finite family of
disjoint sets $A_{t}\in V(t\in T)$ the following implication is
true
\begin{equation}\label{additivity}
    A=\bigcup_T A_{t}\in V\impl \mu(A)=\sum_T\mu(A_{t}).
\end{equation}
Denote by $K(v,Z)$ the space of all vector measures $\mu$ from the
prering $V$ into the space $Z$, such that
\begin{center}
$|\mu(A)|\leq m v(A)\fa A\in V$ and some $m$.
\end{center}
The least constant $m$ satisfying the above inequality is denoted
by $\norm{\mu}$. It is easy to see that the pair $(K(v,Z),
\norm{\mu})$ forms a Banach space.
\end{defin}

\bigskip

Assume that $c_A$ denotes the characteristic function of the set
$A$ that is $c_A(x)=1$ on $A$ and takes value zero elsewhere. Let
$S(V,Y)$ denote the space of all functions of the form
\begin{equation}\label{B}
    h=y_{1}c_{A_1}+\ldots+y_{k}c_{A_{k}},\text{ where }
    y_{i}\in Y, A_{i}\in V.
\end{equation}
The sets $A_i$ in above formula are supposed to be disjoint.
Notice also that we extended the multiplication by scalars by
agreement $y\lambda=\lambda y$ for all vectors $y$ and scalars
$\lambda.$ The family $S(V,Y)$ of functions will be called the
family of {\bf simple functions} generated by the prering $V.$ For
fixed $u\in U$ and $\mu\in K(v,Z)$ define the operator
\begin{center}
    $\int u(h,d\mu)=u(y_{1},\mu(A_{1}))+\ldots+u(y_{k},\mu(A_{k}))$.
\end{center}
Define also
\begin{center}
    $\int hdv=y_{1}v(A_{1})+\ldots+y_{k}v(A_{k})$.
\end{center}
The operators $\int h\,dv$ and $\int u(h,d\mu)$ are well defined,
that is, they do not depend on the representation of the function
$h$ in the form (\ref{B}).

Let $|h|$ denote the function defined by the formula
$|h|(x)=|h(x)|$ for $x\in X$. We see that if $h\in S(V,Y)$, then
$|h|\in S(V,R)$. Therefore the following functional $||h||=
\int|h|d\,v$ is well defined for $h\in S(V,Y)$.

\bigskip

The following development of the theory of Lebesgue and Bochner
summable  functions and of the integrals  are  from Bogdanowicz
\cite{bogdan10}.

\begin{lemma}[Elementary integrals on simple functions]
\label{Elementary integrals} \label{Elementary lemma} \label{Lemma
1} The following statements describe the basic relations between
the notions that we have just introduced.
\begin{enumerate}
\item   The space $S(V,Y)$  is linear, $||h||$
 is a seminorm on it, and $\int h\,dv$  is a
 linear operator on $S(V,Y),$  and $|\int h\,dv|\leq||h||$  for all $h\in S(V,Y)$.
\item   If $g\in S(V,R)$  and $f\in S(V,Y)$,  then $gf\in S(V,Y).$
\item  $\int h dv\geq 0$  if $h\in S(V,R)$  and $h(x)\geq 0$  for
all $x\in X$.

\item   $\int g dv\geq\int fdv$  if $g, f\in S(V,R)$ and $g(x)\geq
f(x)$  for all $x\in X$. \item  The operator $\int u(h,d\mu)$  is
trilinear from the product space $U\times S(V,Y)\times K(v,Z)$
into the space $W$  and
\begin{equation}
    \left|\int u(h,d\mu)\right|\leq|u|\,||h||\,||\mu||
    \qtext{ for all }u\in U,\ h\in S(V,Y),\ \mu\in K(v,Z).
\end{equation}
\end{enumerate}
\end{lemma}

\bigskip


Let $N$ be the family of all sets $A\subset X$ such that for every
$\varepsilon>0$ there exists a countable family $A_{t}\in V(t\in
T)$ such that $A\subset\bigcup_{T}A_{t}$ and
$\sum_{T}v(A_{t})<\varepsilon$. Sets of the family $N$ will be
called {\bf null-sets.} This family represents a {\bf sigma-ideal}
of sets in the power set $\Power(X),$ that is, it has the
following properties: if $A\in N$, then $B\cap A\in N$ for any set
$B\subset X,$ and the union of any countable family of null-sets
$A_{t}\in N(t\in T)$ is also a null-set $\bigcup_T A_{t}\in N$.

A condition $C(x)$ depending on a parameter $x\in X$ is said to be
{\bf satisfied almost everywhere} if there exists a set $A\in N$
such that the condition is satisfied at every
point of the set $X\less A$. 

By a {\bf basic sequence} we shall understand a sequence $s_{n}\in
S(V,Y)$ of functions for which there exists a series with terms $
h_{n}\in S(V,Y)$
 and a constant $M>0$ such that $s_{n}=h_{1}+h_{2}+\ldots+h_{n},$
 where  $||h_{n}||\leq M4^{-n}$ for all
$n=1,2, \ldots$.

\bigskip

\begin{lemma}[Riesz-Egoroff property of a basic sequence]
\label{Riesz-Egoroff property of basic sequence} \label{Basic
sequence lemma} \label{Lemma 2} Assume that $(X,V,v)$ is a
positive measure space on a prering $V$ and $Y$ is a Banach space.
Then the following is true.
\begin{enumerate}
\item {{\rm [Riesz]} \it If} $s_{n}\in S(V,Y)$ {\it is a basic
sequence, then there exists a function} $f$ {\it from} {\it the
set} $X$ {\it into the Banach space} $Y$ {\it and a null-set} $A$
{\it such that} $s_{n}(x)\rightarrow f(x)$ {\it for all} $ x\in$
$X\backslash A$. \item {{\rm [Egoroff]} \it Moreover, for every}
$\varepsilon>0$ {\it and} $\eta>0$, {\it there exists an index}
$k$ {\it and a countable} {\it family of sets} $A_{t}\in V(t\in
T)$ {\it such that}
\begin{center}
$A  \subset\bigcup_{T}A_{t}$ \quad{\it and} \quad $
\sum_{T}v(A_{t})<\eta$
\end{center}
{\it and for every} $n\geq k$
\begin{center}
$|s_{n}(x)-f(x)|<\varepsilon$  \quad{\it if} \quad $x\not\in
\bigcup_{T}A_{t}$.
\end{center}
\end{enumerate}
\end{lemma}

\bigskip

\begin{lemma}[Dunford's Lemma]
\label{Dunford's Lemma} \label{Lemma 3} \label{Basic sequence
converging a.e. to 0 converges in seminorm} Assume that $(X,V,v)$
is a positive measure space on a prering $V$ and $Y$ is a Banach
space. Then the following is true.

If $s_{n}\in S(V,Y)$  is a basic sequence converging almost
everywhere to zero $0$, then the sequence of seminorms $||s_{n}||$
converges to zero.
\end{lemma}

\bigskip

\section{The Spaces of Lebesgue and Bochner Summable Functions}

\bigskip

\begin{defin}[Lebesgue and Bochner spaces]
Assume that $(X,V,v)$ is a measure space over a prering $V$ of an
abstract space $X.$

Let $L(v,Y)$ denote the set of all functions $f:X\mapsto Y,$ such
that there exists basic sequence $s_{n}\in S(V,Y)$ that converges
almost everywhere to the function $f.$

The space $L(v,Y)$ is called the space of {\bf Bochner summable}
functions and, for the case when $Y$ is equal to the space $R$ of
reals, $L(v,R)$ represents the space of {\bf Lebesgue summable}
functions.
\end{defin}

Define
\begin{center}
    $||f||=  \lim_n||s_{n}||,\ \int u(f,d\mu)=
    \lim_n\int u(s_{n},d\mu),\ \int fdv=\lim_n\int s_{n}dv$.
\end{center}
Since the difference of two basic sequences is again a basic
sequence, therefore it follows from the Elementary Lemma
\ref{Lemma 1} and Dunford's Lemma \ref{Lemma 3} that the operators
are well defined, that is, their values do not depend on the
choice of the particular basic sequence convergent to the function
$f$.

\bigskip

\begin{lemma}[Density of simple functions in $L(v,Y)$]
\label{Lemma 4} \label{Density of simple functions in $L(v,Y)$}
Assume that $(X,V,v)$ is a positive measure space on a prering $V$
and $Y$ is a Banach space.
Let  $s_{n}\in S(V,Y)$ be a basic sequence convergent almost
everywhere to a function $f$. Then $||s_{n}-f||\into 0.$
\end{lemma}

\bigskip

\begin{thm}[Basic properties of the space $L(v,Y)$]
\label{Basic properties of the space $L(v,Y)$} \label{Theorem 1}
Assume that $(X,V,v)$ is a positive measure space on a prering $V$
and $Y$ is a Banach space. Then the following is true.

\begin{enumerate}
\item The space $L(v,Y)$ is linear and $||f||$ represents a
seminorm
    being an extension of the seminorm from the space $S(V,Y)$ of simple functions.
\item We have $||f||=0$ if and only if $f(x)=0$ almost everywhere.
\item The functional $||f||$ is $a$
    complete seminorm on $L(v,Y)$ that is given a sequence of functions $f_n\in L(v,Y)$
    such that $\norm{f_n-f_m}\stackrel{nm}\into 0$ there exists a function $f\in L(v,Y)$
    such that $\norm{f_n-f}\stackrel{n}\into 0.$
\item If $f_{1}(x)=f_{2}(x)$ almost everywhere and $f_{2}\in
L(v,Y)$, then $f_{1}\in L(v,Y)$ and
\begin{center}
$||f_{1}||=||f_{2}||,\ \int f_{1}dv=\int f_{2}dv,\ \int
u(f_{1},d\mu)=\int u(f_{2},d\mu)$.
\end{center}
\item The operator $  \int fdv$ is linear and represents an
extension onto $L(v,Y)$ of the operator from $S(V,Y)$. It
satisfies the condition $|\int fdv|\leq||f||$ for all $f\in
L(v,Y)$. \item The operator $  \int u(f,d\mu)$ is trilinear on
$U\times L(v,Y)\times  K(v,Z)$ and represents an extension of the
operator from the space $ U\times S(V,Y)\times K(v,Z).$
 It satisfies the condition:
\begin{center}
$|\int u(f,d\mu)|\leq|u|\,||f||\,||\mu||\fa u\in U, f\in L(v,Y),
\mu\in  K(v,Z).$
\end{center}
\end{enumerate}
\end{thm}


From Theorem \ref{Theorem 1} we see that the obtained integrals
are continuous under the convergence with respect to the seminorm
$\norm{\ },$ that if $||f_{n}-f||\rightarrow 0$, then
\begin{center}
$\int f_{n}dv\rightarrow\int fdv$ and $\int
u(f_{n},d\mu)\rightarrow\int u(f,d\mu)$.
\end{center}
The following theorem characterizes convergence with respect to
this seminorm.


\begin{thm}[Characterization of the seminorm convergence]
\label{Characterization of the seminorm convergence}
\label{Theorem 2} Assume that $(X,V,v)$ is a positive measure
space on a prering $V$ and $Y$ is a Banach space.

Assume that we have a sequence of summable functions $f_{n}\in
L(v,Y)$ and some function $f$ from the set $X$ into the Banach
space $Y.$ Then the following conditions are equivalent
\begin{enumerate}
\item The sequence $f_n$ is Cauchy, that is
$||f_{n}-f_{m}||\stackrel{nm}\into 0,$
 and there exists a subsequence $f_{k_{n}}$  convergent
 almost everywhere to the function $f$.
\item The function $f$ belongs to the space $L(v,Y)$  and
$||f_{n}-f||\rightarrow 0.$
\end{enumerate}
\end{thm}

When the space $Y=R,$ the space $L(v,R)$ represents the space of
Lebesgue summable functions. We have the following relation
between Bochner summable functions and Lebesgue summable
functions.

\bigskip

\begin{thm}[Norm of Bochner summable function is Lebesgue summable]
\label{Norm of Bochner summable function is Lebesgue summable}
\label{Theorem 3} Let $(X,V,v)$ be a positive measure space on a
prering $V$ and assume that $Y$ is a Banach space.

If $f$ belongs the the space $L(v,Y)$ of Bochner summable
functions, then the function $|f|=|f(\cdot)|$ belongs to the space
$L(v,R)$ of Lebesgue summable functions and we have the identity
\begin{center}
    $||f||=\int|f(\cdot)|\,dv\fa f\in L(v,Y).$
\end{center}
\end{thm}


\begin{thm}[Properties of Lebesgue summable functions]
\label{Properties of Lebesgue summable functions} \label{Theorem
4} Let $(X,V,v)$ be a positive measure space on a prering $V$ and
$L(v,R)$ the Lebesgue space of $v$-summable functions.
\begin{description}
    \item[(a)]  If $f\in L(v,R)$  and
        $f(x)\geq 0$ almost everywhere on $X$ then $\int fdv\geq 0$.
    \item[(b)]  If $f,g\in L(v,R)$
        and $f(x)\geq g(x)$ almost everywhere on $X$ then $\int f\,dv\geq\int g\,dv.$
    \item[(c)]  If $f,g\in L(v,R)$ and $h(x)=\sup\{f(x),g(x)\}$ for all $x\in X$
        then $h\in L(v,R).$
    \item[(d)]  Let $f_{n}\in L(v,R)$ be a monotone sequence with respect
        to the relation less or equal almost everywhere.
        Then there exists a function $f\in L(v,R)$
        such that $f_{n}(x)\into f(x)$ almost everywhere on $X$
        and $||f_{n}-f||\rightarrow 0$ if and only if the sequence of numbers
        $\int f_{n}\,dv$ is bounded.
    \item[(e)] Let $g,f_{n}\in L(v,R)$ and $f_{n}(x)\leq g(x)$ almost
        everywhere on $X$ for $n=1,2, \ldots$.
        Then the function $h(x)=\sup\{f_{n}(x):\ n=1,2,\ldots\}$
        is well defined almost everywhere on $X$
        and is summable, that is,  $h\in L(v,R)$.
        A function defined almost everywhere is said to be
        summable if it has a summable extension onto the space $X.$
\end{description}
\end{thm}

From part (d) of the above theorem we can get the classical
theorem due to Beppo Levi \cite{levi}.

\begin{cor}[Beppo Levi's Monotone Convergence Theorem]
\label{Beppo Levi's Monotone Convergence Theorem} \label{Monotone
Convergence Theorem} Assume that $(X,V,v)$ is a positive measure
space on a prering $V$ and $L(v,R)$ the Lebesgue space of
$v$-summable functions.

Let $f_{n}\in L(v,R)$ be a monotone sequence with respect to the
relation less or equal almost everywhere. Then there exists a
function $f\in L(v,R)$ such that $f_{n}(x)\into f(x)$ almost
everywhere on $X$ and $$\int f_{n}dv\rightarrow \int f\,dv$$ if
and only if the sequence of numbers $\int f_{n}\,dv$ is bounded.
\end{cor}

\begin{thm}[Lebesgue's Dominated Convergence Theorem]
\label{Dominated Convergence Theorem} \label{Theorem 5} Let
$(X,V,v)$ be a positive measure space on a prering $V$ and $Y$ a
Banach space.

Assume that we are given a sequence $f_{n}\in L(v,Y)$ of Bochner
summable functions that can be majorized by a Lebesgue summable
function $g\in L(v,R),$ that is for some null set $A\in N$ we have
the estimate
$$|f_{n}(x)|\leq g(x)\fa x\not\in A\text{ and } n=1,2, \ldots$$
Then the condition
   $$f_{n}(x)\rightarrow f(x)\qtext{a.e. on X}$$
implies the relations
    $$f\in L(v,Y)\qtext{and} ||f_{n}-f||\rightarrow 0$$
and, therefore, also the relations
\begin{center}
    $\int f_n\,dv\into \int f\,dv\qtext{and}\int u(f_n,d\mu)\into \int u(f,d\mu)$
\end{center}
for every bilinear continuous operator $\ u\ $ from the product
$Y\times Z$ into the Banach space $W$ and any vector measure
$\mu\in K(v,Z).$
\end{thm}


\section{Continuous functions on compact sets are summable}

\bigskip

Now let $(X,V,v)$ denote one of the following  measure spaces:
$(R,V,v)$ be the Riemann measure space, or the space $R^n$ with the prering
consisting of all the cubes of the form
\begin{equation*}
 A=J_1\times\ldots\times J_n,\ J_i=(a_i,b_i],\ a_i\le b_i
\end{equation*}
and Riemann measure
\begin{equation*}
 v(A)=(b_1-a_1)\cdots(b_n-a_n),
\end{equation*}
or let $X$ be a topological Hausdorff space, the prering $V$
consists of sets of the form $A=Q_1\less Q_2$ where $Q_i$
are compact sets, and the measure $v$ be any countably additive
nonnegative finite-valued function on $V.$

In the case when $X$ is a locally compact topological
group the Haar measure restricted to the prering $V$
provides a nontrivial example of such measure space.
For details see Halmos \cite{halmos}, Chapters 10 and 11.

The following is Theorem 8, page 498, of Bogdanowicz \cite{bogdan10}.

\begin{thm}[Summability of continuous functions on compact sets]
Assume that the triple $(X,V,v)$ represents
one of the above measure spaces and $f$ a continuous
function from a compact set $K\subset X$ into a Banach space $Y.$
Then the function $f$ is Bochner summable on $K$ that is we have
$c_Kf\in L(v,Y)$ and thus the integral $
    \int_Kf(t)\,dt
$ exists.
\end{thm}

\bigskip

\bp
    Consider first the case of the space $R^n.$
    The set $K$ being compact is bounded. Thus
    there exists a cube $I$ containing the set $K.$
    Divide the cube $I$ into finite number of
    disjoint cubes of diameter less  than
    $\frac{1}{n}.$ Let $I^n_j\ (j=1,2,\ldots,k_n)$ be the cubes such that
    $K\subset \bigcup_jI^n_j$ and $K\cap I^n_j\neq \emptyset.$

    In the case of a topological space $X$ the set $f(K),$ as
    an image of a compact set by means of a continuous function,
    is compact. There exist disjoint sets  $I^n_j\ (j=1,2,\ldots,k_n)$
    of diameter $<\frac{1}{n}$ such that $\bigcup_jI^n_j=f(K)$
    and each of the sets $I^n_j$ is the intersection of an open
    set with a closed set. Therefore $f^{-1}(I^n_j)\in V.$

    In any case choose $x^n_j\in K\cap I^n_j.$ and let
    \bb
        s_n=\sum_jf(x^n_j)c_{I^n_j}\fa n=1,2,\ldots
    \ee
    The sequence $s_n$ consists of simple functions thus $s_n\in L(v,Y).$
    It converges everywhere on $X$ to the function $c_Kf$ and is dominated
    everywhere by the simple function $m\,c_I$  or $m\,c_K,$ where
    \begin{equation*}
         m=\sup\set{\abs{f(x)}:\,x\in K}.
    \end{equation*}
    Hence from the Lebesgue Dominated Convergence Theorem we get $c_Kf\in L(v,Y).$
\ep

The above theorems are from Bogdanowicz \cite{bogdan10} in the
order as they have been proved in that paper.

\bigskip

\section{Summable sets form a delta ring}

\bigskip

Assume now again that we have a measure space $(X,V,v)$ on a
prering $V$ of subsets of an abstract space $X.$ Following
Bogdanowicz \cite{bogdan15} and \cite{bogdan20} denote by $V_c$
the family of all sets $A\subset X$ whose characteristic function
$c_A$ is $v$-summable that is $c_A\in L(v,R).$ Put $v_c(A)=\int
c_Adv$ for all sets $A\in V_c.$ From the properties of the
Lebesgue summable functions (\ref{Properties of Lebesgue summable
functions}) we can deduce the following proposition.


\bigskip

\begin{pro}[Summable sets form a delta ring]
Assume that $(X,V,v)$ is  a measure space on a prering $V$ of
subsets of an abstract space $X.$

Then the family $V_c$ of summable sets forms a $\delta$-ring and
the set function $v_c$ forms a measure. If in addition $X\in V_c$ then $V_c$ forms
a $\sigma$-algebra.
\end{pro}


\bp
    First of all notice that from Theorem (\ref{Theorem 3}) follows that
    absolute value of a Lebesgue summable function is itself summable
    \begin{equation}\label{summability of absolute value}
    |f|\in L(v,R)\fa f\in L(v,R).
    \end{equation}
    Thus by linearity of the space $L(v,R)$ of Lebesgue summable functions
    and from the identities
    \begin{equation*}
    \begin{split}
    (f\vee g)(x)=&\sup\set{f(x),g(x)}=\frac{1}{2}(f(x)+g(x)+|f(x)-g(x)|)\qtext{and}\\
    (f\wedge g)(x)=&\inf\set{f(x),g(x)}=\frac{1}{2}(f(x)+g(x)-|f(x)-g(x)|)
    \end{split}
    \end{equation*}
    follows that $f\vee g$ and $f\wedge g$ are in $L(v,R)$ if $f,g$ are in $L(v,R).$
    Thus from the identities
    \begin{equation*}
    c_{A\cup B}=c_A\vee c_B\qtext{and}c_{A\cap B}=c_A\wedge c_B\qtext{and}c_{A\less B}=
    c_A-c_A\wedge c_B
    \end{equation*}
    we can conclude that the family $V_c$ of summable sets forms a ring.

    Now if $A_n\in V_c$ is a sequence of summable sets
    and $B_n=\bigcap_{j\le n} A_j$ and $B=\bigcap_{j\ge 1} A_j,$
    from the Dominated Convergence Theorem (\ref{Dominated Convergence Theorem})
    and from the relations
    \begin{equation*}
    |c_{B_n}(x)|=c_{B_n}(x)\le c_{A_1}(x)\qtext{and}c_{B_n}(x)\into c_B(x)\fa x\in X
    \end{equation*}
    we get that $B\in V_c.$ Thus the family $V_c$ of summable sets forms a
    $\delta$-ring.

    In the case, when $X\in V,$ we get from the de Morgan law
    and the fact that $V_c$ forms a $\delta$-ring that
    \begin{equation*}
    \bigcup_{n\ge 1}A_n=X\less \bigcap_{n\ge 1}(X\less A_n)\in V_c.
    \end{equation*}
    Hence in this case $V_c$ forms a $\sigma$-algebra.

    To show that the triple $(X,V_c,v_c)$ forms a positive
    measure space assume that $A\in V_c$
    and a sequence of disjoint sets $A_n\in V_c$ forms a decomposition of the set $A.$
    So $A=\bigcup_{j\ge 1} A_j.$ Let $B_n=\bigcup_{j\le n} A_j.$
    From the Dominated Convergence Theorem (\ref{Dominated Convergence Theorem})
    and from the relations
    \begin{equation*}
    |c_{B_n}(x)|=c_{B_n}(x)\le c_{A}(x)\qtext{and}c_{B_n}(x)\into c_A(x)\fa x\in X
    \end{equation*}
    and linearity of the integral,
    we get that
    \begin{equation*}
    \begin{split}
    v_c(A)&=\lim_n v_c(B_n)=\lim_n \int c_{B_n}dv\\&=\lim_n \sum_{j\le n}\int c_{A_j}dv=
    \lim_n \sum_{j\le n}v_c(A_j)=\sum_{j}v_c(A_j).
    \end{split}
    \end{equation*}
     Thus $v_c$ is countably additive on the delta ring $V_c.$
\ep


Now let us consider the case when the space $Y$ is the space $R$
of reals and $Z$ any Banach space and the bilinear operator $u$ is
the multiplication operator $u(r,z)=rz.$

\bigskip

\begin{pro}[Isomorphism and isometry of $K(v,Z)$ and $K(v_c,Z)$]
Assume that $(X,V,v)$ is  a measure space on a prering $V$ of
subsets of an abstract space $X.$

Every vector measure $\mu\in K(v,Z)$ can be extended from the
prering $V$ onto the delta ring $V_c$ by the formula
\begin{equation*}
    \mu_c(A)=\int u(c_A,d\mu)\fa A\in V_c.
\end{equation*}
This extension establishes isometry and isomorphism between the
Banach spaces $K(v,Z)$ and $K(v_c,Z).$
\end{pro}

\bigskip

\section{A characterization of summable functions}

\bigskip

\begin{pro}[A characterization of Bochner summable functions]
\label{A characterization of Bochner summable functions} Assume
that $(X,V,v)$ is  a measure space on a prering $V$ of subsets of
an abstract space $X $ and $Y$ a Banach space.

A function $f$ from $X$ into $Y$ belongs to the space $L(v,Y)$ if
and only if there exist a sequence $s_n\in S(V,Y)$ of simple
functions and a non-negative summable function $g\in L(v,R)$ such
that $s_n(x)\into f(x)$ almost everywhere on $X$ and
\begin{equation*}
    |s_n(x)|\le g(x)\fa n=1,2,\ldots \qtext{and}\qtext{almost all} x\in X.
\end{equation*}
\end{pro}

\bigskip

\bp
    If $f\in L(v,Y)$ then there exists a basic sequence of the form
    \bb
        s_n=h_1+h_2+\cdots+h_n
    \ee
    converging almost everywhere to the function $f.$
    Notice that the sequence
    \bb
        S_n=|h_1|+|h_2|+\cdots+|h_n|
    \ee
     is nondecreasing and is basic. Thus it converges almost
     everywhere to some summable function $g\in L(v,R).$
    Since
    \bb
        |s_n(x)|\le S_n(x)\le g(x)\qtext{for almost all}x\in X,
    \ee
    we get the necessity of the condition.

    The sufficiency of the condition follows from the Dominated Convergence Theorem.
\ep

\bigskip

\begin{pro}[Summability of a product of functions]
\label{Summability of a product of functions} Assume that
$(X,V,v)$ is  a measure space on a prering $V$ of subsets of an
abstract space $X $ and $\ Y$ a Banach space. Let
\begin{equation*}
    u(y)=\frac{1}{|y|}\,y\text{ if }|y|>0\qtext{and}u(y)=0\text{ if }|y|=0.
\end{equation*}
Assume that $f\in L(v,Y)$ and $g\in L(v,R).$ Then the product
function $u\circ f\cdot g$ is summable that is $u\circ f\cdot g\in
L(v,Y),$ where $u\circ f$ denotes the composition $(u\circ
f)(x)=u(f(x))$ for all $x\in X.$
\end{pro}

\bigskip

\bp
    Take any natural number $k$ and define a function
    \bb
        u_k(y)=(k|y|\wedge 1)\frac{1}{|y|}y\fa y\in Y,\,|y|>0\qtext{and}u_k(0)=0.
    \ee
    Notice that the function $u_k$ is continuous and
    \bb
        \lim_k u_k(y)=u(y) \qtext{ and } |u_k(y)|\le 1\fa y\in Y.
    \ee
    Let $s_n$ be a basic sequence converging almost everywhere to $f$ and
    $S_n$ a basic sequence converging almost everywhere to $g.$
    Let $G\in L(v,R)$ be a majorant for the sequence $S_n.$
    Then we have that the sequence
    \bb
        h_{k\,n}=u_k\circ s_n\cdot S_n\in S(V,Y)
    \ee
    consists of simple functions and when $n\into \infty$
    it converges almost everywhere to the function
    \bb
        h_k=u_k\circ f\cdot g.
    \ee
    Since $G$ majorizes the sequence $h_{k\,n},$ from the Dominated Convergence
    Theorem we get $h_k\in L(v,Y)$ and moreover
    \bb
        |h_k(x)|\le G(x)\qtext{for almost all}x\in X.
    \ee
    Passing to the limit $k\into\infty$ and applying the Dominated Convergence
    Theorem yields that $$u\circ f\cdot g\in L(v,Y).$$
\ep



\begin{defin}[Summability on sets]
Assume that $(X,V,v)$ is  a measure space on a prering $V$ of
subsets of an abstract space $X $ and $Y$ a Banach space.

We shall say that a function $f:X\mapsto Y$ is {\bf summable on a
set} $A\subset X$ if the product function $c_Af$ is summable and
we shall write
$$\int_Af\,dv=\int c_Af\,dv.$$
\end{defin}

\begin{defin}[Vector measure of finite variation]
A vector measure $\mu$ from a prering $V$ into a Banach space $Y$
is said to be of {\bf finite variation} on $V$
if 
\bb
    |\mu|(A)=\sup\sum_{t\in T}|\mu(A_t)|<\infty\fa A\in V
\ee where the supremum is taken over all finite disjoint
decompositions $A_t\in V\,(t\in T)$ of the set $A=\bigcup_{t\in
T}A_t.$ The set function $|\mu|$ is called the {\bf variation} of
the vector measure $\mu.$
\end{defin}

\bigskip


\begin{pro}[Sets on which a function is summable form a $\delta$-ring]
Assume that $(X,V,v)$ is  a measure space on a prering $V$ of
subsets of an abstract space $X $ and $Y$ a Banach space and
assume that $f$ is an arbitrary function from $X$ into $Y.$ Denote
by $V_f$ the family of all sets $A\subset X$ on which the function
$f$  is summable.

If $f:X\mapsto Y$ is an arbitrary function then $V_f$ forms a
$\delta$-ring and $\mu(A)=\int_Af\,dv$ forms a $\sigma$-additive
vector measure of finite variation on $V_f.$
\end{pro}

\bigskip

\bp
    Assume that $A,B\in V_f.$ Then $c_Af\in L(v,Y)$ and $|c_Bf|\in L(v,R).$
    From Proposition \ref{Summability of a product of functions} we get
    \bb
        c_{A\cap B}f=c_Ac_B\,u\circ{f}\cdot|f|=u\circ (c_Af)\cdot |c_Bf|\in L(v,Y).
    \ee
    Thus $A\cap B\in V_f.$ It follow from linearity of the space $L(v,Y)$
    and the identities
    \bb
        c_{A\less B}=c_A-c_{A\cap B}\qtext{and}
        c_{A\cup B}=c_{A\less B}+c_{A\cap B}+c_{B\less A}
    \ee
    that $V_f$ forms a ring.

    Now using the Dominated Convergence Theorem we can easily prove that
    $V_f$ forms a $\delta$-ring and the set function
    \bb
        \mu(A)=\int_Af\,dv\fa A\in V_f
    \ee
    forms a $\sigma$-additive vector measure and
    \bb
        |\mu|(A)\le \int_A|f|\,dv<\infty\fa A\in V_f.
    \ee
\ep

\bigskip

The family of sets $V_f$ on which a function $f$ is summable may
consist only of the empty set. However in the case of a summable
function this family is rich as follows from the following
corollary.

\bigskip


\begin{cor}[The collection of sets on which a summable
function is summable forms a  $\sigma$-algebra]
Assume that $(X,V,v)$ is  a measure space on a prering $V$ of
subsets of an abstract space $X $ and $Y$ a Banach space.

If $f\in L(v,Y)$ is a summable function then the family $V_f$ of
sets, on which $f$ is summable, forms a $\sigma$-algebra
containing all summable sets that is we have the inclusion
$$V\subset V_c\subset V_f$$
and the set function $\mu(A)=\int_Af\,dv$ is $\sigma$-additive of
finite total variation $|\mu|(X)\le \int|f|\,dv.$
\end{cor}

\bigskip

\bp
    To prove this corollary notice that similarly as before we can
    prove that product $gf$ of a summable bounded function $g\in L(v,R)$
    with a summable function $f\in L(v,Y)$ is summable $gf\in L(v,Y).$
    This implies that $V\subset V_c\subset V_f.$
\ep

\bigskip

For further studies of vector measures we recommend Dunford and
Schwartz \cite{DS1}, and for extensive survey of the state of the
art in the theory of vector measures see the monograph of Diestel
and Uhl \cite{diestel-uhl}.

\bigskip

\section{Extensions to Lebesgue measures}

\bigskip

If $V$ is any nonempty collection of subsets of an abstract space
$X$ denote by $V^{\sigma}$ the collection of sets that are
countable unions of sets from $V$ and denote by $V^{r}$ the
collection \bb
    V^r=\set{A\subset X:\ A\cap B\in V\fa B\in V}.
\ee Now assume that $(X,V,v)$ is a measure space and $V_c$ the
$\delta$-ring of summable sets and $v_c(A)=\int c_Adv.$ It is easy
to prove that $V_c^\sigma$ forms the smallest $\sigma$-ring
containing the prering $V$ and the family $N$ of $v$-null sets.
Moreover the set function defined by
\begin{equation}\label{completing v_c to Lebesgue measure}
    \mu(A)=\sup\set{v_c(B):\ B\subset A,\,B\in V_c}\fa A\in V_c^\sigma
\end{equation}forms a Lebesgue measure on $V_c^\sigma.$

A Lebesgue measure is called {\bf complete} if all subsets of sets
of measure zero are in the domain of the measure and thus have
measure zero. The above Lebesgue measure $\mu$ can be
characterized as the smallest extension of the measure $v$ to a
complete Lebesgue measure. Hence this measure is unique.

One can prove that \bb
    V_c^r=\bigcap_{f\in L(v,Y)}V_f.
\ee The family $V_c^r$ as an intersection of $\sigma$-algebras
forms itself a $\sigma$-algebra containing the $\delta$-ring
$V_c.$ The smallest $\sigma$-algebra $V^a$ containing $V_c$ is
given by the formula \bb
    V^a=\set{A\subset X:\ A\in V_c^\sigma\text{ or }X\less A\in V_c^\sigma}.
\ee

If $X\in V_c^\sigma$ then the sigma algebras coincide
$V_c^\sigma=V^a=V_c^r.$ If  $X\not\in V_c^\sigma,$ one can always
extend the measure $v$ to a Lebesgue measure on  $V^a$ or $V_c^r$
by the formula
\begin{equation}\label{trivial extension}
    \mu(A)=v_c(A)\text{ if } A\in V_c
    \qtext{and}\mu(A)=\infty\text{ if }A\not\in V_c.
\end{equation}
However if $\sup\set{v_c(A):\ A\in V_c}=a<\infty$ and $X\not\in
V_c^\sigma$ the extensions are not unique. Indeed one can take in
this case $\mu(X)=b,$ where $b$ is any number from the interval
$\langle a,\infty),$ and put \bb
    \mu(A)=v_c(A)\text{ if }A\in V_c\qtext{and}
    \mu(A)=b-v_c(X\less A)\text{ if }X\less A\in V_c
\ee to extend  the measure $v_c$ onto the $\sigma$-algebra $V^a$
preserving sigma additivity.

Consider an example. Let $(X,V,v)$ be the following measure space:
\bb
    X=R,\quad V=\set{\emptyset,\{n\}:\ n=1,2,\ldots},\quad v(\emptyset)=0,\,v(\{n\})=2^{-n}.
\ee In this case the family $N$ of null sets contains only the
empty set $\emptyset,$ the family $S$ of simple sets consists of
finite subsets of the set of natural numbers $\N,$ the family
$V_c$ of summable sets consists of all subsets of $\N,$ we have
$V_c^\sigma=V_c,$ the smallest $\sigma$-algebra extending $V_c$
consists of sets that either are subsets of $\N$ or their
complements are subsets of $\N,$ finally $V_c^r=P(R)$ consists of
all subsets of $R.$ Since $v_c(\N)=1$ is the supremum of $v_c,$
the measure $v_c$ has many extensions onto the $\sigma$-algebras
$V^a$ and $V_c^\sigma.$ One extension onto $P(R)$ is given by
\ref{trivial extension} and another, for instance, by \bb
    \mu(A)=\sum_{n\in A\cap \N}2^{-n}\fa A\subset R.
\ee

\bigskip

In view of the existence of a variety of extensions of a measure
from a prering onto $\delta$-rings and the multiplicity of
extensions to Lebesgue measures it is important to be able to
identify measures that generate the same class of Lebesgue-Bochner
summable functions $L(v,Y)$ and the same trilinear integral $\int
u(f,\,d\mu)$ and thus the ordinary Bochner integral $\int f\,dv.$
In this regard we have the following theorems.

Assume that $(X,V_j,v_j),\,(j=1,2)$ are two measure spaces over
the same abstract space $X$ and $Y,Z,W$ are any Banach spaces and
$U$ is the Banach space of bilinear bounded operators from the
product $Y\times Z$ into $W.$

\bigskip

\begin{thm}[When $L(v_2,Y)$ extends $L(v_1,Y)$?]
\label{When $L(v_2,Y)$ extends $L(v_1,Y)$?} For every Banach space
$Y$ we have $L(v_1,Y)\subset L(v_2,Y)$ and \bb
    \int f\,dv_1=\int f\,dv_2\fa f\in L(v_1,Y)
\ee if and only if $V_{1c}\subset V_{2c}$ and \bb
    v_{1c}(A)=v_{2c}(A)\fa A\in V_{1c}
\ee that is the measure $v_{2c}$ represents an extension of the
measure $v_{1c}.$
\end{thm}

\bigskip

Consequently we have the following theorem.

\begin{thm}
For any Banach space $Y$ and any bilinear bounded transformation
$u\in U$ we have $L(v_1,Y)= L(v_2,Y)$ and \bb
    \int f\,dv_1=\int f\,dv_2\fa f\in L(v_1,Y)
\ee and the spaces of vector measures
$K(v_1,Z),K(v_2,Z),K(v_{1c},Z),K(v_{2c},Z)$ are isometric and
isomorphic and \bb
    \int u(f,d\mu_1)=\int u(f,d\mu_2)=\int u(f,d\mu_{1c})
    =\int u(f,d\mu_{2c})\fa f\in L(v_1,Y),
\ee where $\mu_1,\mu_2,\mu_{1c},\mu_{2c}$ are vector measures that
correspond to each other through the isomorphism, if and only if,
the completions of the measures $v_1,v_2$ coincide
$v_{1c}=v_{2c}.$
\end{thm}

For proofs of the above theorems see Bogdanowicz \cite{bogdan20}.
It is important to relate the above theorems to the classical
spaces of Lebesgue and Bochner summable functions and the
integrals generated by Lebesgue measures. Since there is a great
variety of approaches to construct these spaces we shall
understand by a classical construction of the Lebesgue space
$L(\mu,R)$ the construction developed in Halmos \cite{halmos} and
by classical approach to the theory of the space $L(\mu,Y)$ of
Bochner summable functions as presented in Dunford and Schwartz
\cite{DS1}.

Now if $(X,V,v)$ is a measure space on a prering $V$ and
$(X,M,\mu)$ represents a Lebesgue measure space where $\mu$ is the
smallest extension of the measure $v$ to a Lebesgue complete
measure on the $\sigma$-ring $M,$ then we have the following
theorem.

\begin{thm}
For every Banach space $Y$ the spaces $L(v,Y)$ and $L(\mu,Y)$
coincide and we have \bb
    \int_A f\,dv=\int_A f\,d\mu \fa f\in L(\mu,Y) \qtext{and} A\in M.
\ee
\end{thm}

The above theorem is a consequence of the theorems developed in
Bogdanowicz \cite{bogdan14}.

\bigskip

\section{Tensor product of measure spaces}

\bigskip

Assume now that we have two measure spaces $(X_i,V_i,v_i)$ over
abstract spaces $X_i$ for $i=1,2.$ Consider the Cartesian product
$X_1\times X_2.$ By {\bf tensor product} $V_1\otimes V_2$ of the
families $V_i$ we shall understand the family of sets \bb
    V_1\otimes V_2=\set{A_1\times A_2:\ A_1\in V_1\text{ and }A_2\in V_2}.
\ee We shall use a shorthand notation
    $$v_1\otimes v_2(A_1\times A_2)=v_1(A_1)v_2(A_2)\fa A_i\in V_i\,(i=1,2).$$
\begin{thm}[Tensor product of measure spaces]
Assume that $(X_i,V_i,v_i)$ are measures over abstract spaces
$X_i$ and $V_i$ are prerings for $i=1,2.$ Let the triple $(X,V,v)$
consist of $X=X_1\times X_2,$ $V=V_1\otimes V_2,$ and
$v(A)=v_1(A_1)v_2(A_2)$ for all $A=A_1\times A_2\in V.$ Then $V$
forms a prering and $v$ a $\sigma$-additive finite-valued positive
measure, that is the triple $(X,V,v)$ forms a measure space.
\end{thm}

\bigskip
\bp
    Notice that the following two properties of a family $V$ of subsets
    of a space $X$ are equivalent:
    \begin{itemize}
        \item The family $V$ forms a prering.
        \item The empty set belongs to $V$ and for every two sets $A,B\in V$
                there exists a finite disjoint refinement from the family $V,$
                that is, there exists a finite collection $\{D_1,\ldots,D_k\}$
                of disjoint sets from $V$ such that each of the two sets $A,B$
                can be represented as a union of some sets from the collection.
    \end{itemize}

    Clearly the tensor product $V_1\otimes V_2$ of the prerings contains the empty set.
    Now take any pair of sets $A,B\in V_1\otimes V_2.$ We have $A=A_1\times A_2$
    and $B=B_1\times B_2.$ If one of the sets  $A_1,A_2,B_1,B_2$ is empty then the
    pair $A,B$ forms its own refinement from $V.$ So consider the case when
    all the sets $A_1,A_2,B_1,B_2$ are nonempty.

    Let $C=\set{C_j\in V_1:\ j\in J}$ be a refinement of the
    pair $A_1,B_1$ and $$D=\set{D_k\in V_2:\ k\in K}$$ a refinement of $A_2,B_2.$
    We may assume that the refinements do not contain the empty set.

    The collection of sets $C\otimes D$
    forms a refinement of the pair $A,B.$
    Indeed each set of the pair $A_1,B_1$ can be uniquely
    represented as the union of sets
    from the refinement $C.$ Similarly each set of the pair $A_2,B_2$ can be represented
    in a unique way as union of sets from the refinement $D.$ Since
    the sets of the collection  $C\otimes D$ are disjoint and nonempty,
    each set of the pair
    $A_1\times A_2$ and $B_1\times B_2$ can be uniquely represented as the union
    of the sets from  $C\otimes D.$
    Thus $V=V_1\otimes V_2$ is a prering.

    To prove  that the set function $v=v_1\otimes v_2$
    is $\sigma$-additive take any set $A\times B$
    in $V$ and let $A_n\times B_n\in V$ denote a sequence of disjoint sets whose
    union is the set $A\times B.$ Notice the identity

    \begin{equation}\label{tensor decomposition}
        c_A(x_1)\,c_B(x_2)=\sum_n c_{A_n}(x_1)\,c_{B_n}(x_2)\fa x_1\in X_1,x_2\in X_2.
    \end{equation}
    Fixing $x_2$ and integrating with respect to $v_1$
    both sides of the equation (\ref{tensor decomposition} ) on the basis
    of the Monotone Convergence Theorem, we get
    \begin{equation*}
        v_1(A)\,c_B(x_2)=\sum_n v_1({A_n})\,c_{B_n}(x_2)\fa x_2\in X_2.
    \end{equation*}
    Integrating the above term by term with respect to $v_2$
    and applying again the Monotone Convergence Theorem yields
    \begin{equation*}
        v_1(A)\,v_2(B)=\sum_n v_1({A_n})\,v_2({B_n})
    \end{equation*}
    that is the set function
        $$v(A\times B)=v_1\otimes v_2(A\times B)
        =v_1(A)v_2(B)\fa A\times B\in V_1\otimes V_2$$
    is $\sigma$-additive.
    Hence the triple
        $$(X,V,v)=(X_1\times X_2,V_1\otimes V_2,v_1\otimes v_2)$$
    forms a measure space.
\ep

\bigskip

The above theorem has an immediate generalization to any finite
number of measure spaces.

\bigskip

\begin{thm}[Tensor product of $n$ measure spaces]
Assume that $(X_i,V_i,v_i)$ are measures over abstract spaces
$X_i$ and $V_i$ are prerings for $i=1,\ldots,n.$ Let the triple
$(X,V,v)$  consist of $X=X_1\times\cdots\times X_n,$
$V=V_1\otimes\cdots\otimes V_n,$ and
    $$v(A)=v_1\otimes\cdots\otimes v_n(A)=v_1(A_1)\cdots v_n(A_n)$$
for all $A=A_1\times\cdots\times A_n\in V.$ Then the triple
$(X,V,v)$ forms a measure space.
\end{thm}

\bigskip


\begin{defin}[Classical Lebesgue measure over $R^n$]
To construct the classical Lebesgue measure $\mu$ over the space $R^n,$ 
first take the tensor product of $n$ copies of the Riemann measure
space $(R,V,v),$ and complete it by means of the formula
(\ref{completing v_c to Lebesgue measure}).
\end{defin}

\bigskip

\section{Integration over the space $R$ of reals}

\bigskip

Now let $X=I$ be a closed bounded interval, and let $V$ denote the
prering of all subintervals of $I,$ and $v(A)$ the length of the
interval $A\subset I.$ Clearly the space $(X,V,v)$ is a measure
space as a subspace of the Riemann measure space.


For the case of Riemann measure space we shall use the customary
notation for the integral of a Bochner summable function $f\in
L(v,Y).$ We shall write
\begin{equation*}
    \begin{split}
    \int_{t_1}^{t_2}f(t)\,dt&=\quad\int c_{[t_1,t_2]}f\,dv
        \qtext{if}t_1\le t_2,\ t_1,t_2\in I,\\
    \int_{t_1}^{t_2}f(t)\,dt&=-\int c_{[t_2,t_1]}f\,dv
        \qtext{if}t_1> t_2,\ t_1,t_2\in I.\\
    \end{split}
\end{equation*}

\bigskip


Adopting the above notation yields a convenient formula for any
$f\in L(v,Y)$
\begin{equation*}
    \int_{t_1}^{t_2}f(t)\,dt+\int_{t_2}^{t_3}f(t)\,dt+\int_{t_3}^{t_1}f(t)\,dt=0\fa
    t_1,t_2,t_3\in I.
\end{equation*}

\bigskip

\section{Right-sided antiderivatives and fundamental theorem of calculus}

\bigskip

We remind the reader that the notion of a set of Lebesgue measure zero is the same
as the notion of the null set corresponding to the Riemann measure space over
the reals.

We also use the notion of locally bounded, or locally essentially bounded, or locally
summable over an open set $J$ to mean that each point of the set $J$ has a neighborhood
on which this property holds. Notice that the notion that a function $f$ is locally
bounded, or locally essentially bounded, locally summable on $J$ is equivalent to
the property that the function is bounded, or essentially bounded, or summable
on every compact subset $F$ of the set $J.$

\begin{thm}[If $f'_r$ is locally essentially bounded then fundamental theorem of calculus holds]
Let $J$ be an open interval and $\ Y$ a Banach space.
Assume that $f,g$ are defined on $J$ into $Y.$
If the function $f$ represents a right-sided antiderivative of $g$
over the interval $J$ and $g$ is locally bounded on $J\less T$
where $T$ is a set of Lebesgue measure zero,
then $g$ is Bochner summable and we have
\begin{equation*}
 \int_{t_1}^{t_2}g(u)\,du=f(t_2)-f(t_1)\fa t_1,t_2\in J.
\end{equation*}
\end{thm}

\bigskip

\bp
    Take any two points $t_1,t_2\in J.$ We may assume without loss of generality
    that $t_1<t_2.$ Select $\delta>0$ so that $t_2+\delta\in J.$

    Take any sequence $h_n\in (0,\delta)$ such that $h_n\into 0$ and
    consider the functions
    \begin{equation*}
        g_n(t)=(h_n)^{-1}[f(t+h_n)-f(t)]\fa t\in I=[t_1,t_2].
    \end{equation*}
    Notice that the functions $g_n$ are continuous on the interval $I$
    and since $I$ is compact the functions $g_n$ are summable on $I.$

    Since neglecting a set $T$ of Lebesgue measure zero the function $g$
    is local bounded on $J,$ from compactness of the interval $I$
    we get that there is a constant $m$ such that
    \begin{equation*}
     \norm{g(t)}\le m\fa t\in [t_1,t_2+\delta]\less T.
    \end{equation*}
    Thus from the Strong Mean-Value theorem (\ref{str mean value} follows that
    \begin{equation*}
     \norm{g_n(t)}\le m\fa t\in [t_1,t_2+\delta]\less T.
    \end{equation*}
    Since by definition of right-sided antiderivative
    \begin{equation*}
        g_n(t)\into f'_r(t)\fa t\in I\qtext{ and }f'_r(t)=g(t)\fa t\in I\less T
    \end{equation*}
    from the Dominated Convergence theorem (\ref{Dominated Convergence Theorem})
    we get that the function $g$ is Bochner summable on the interval $I$
    and we have the convergence
    \begin{equation*}
    \int_{t_1}^{t_2}g_n(u)\,du\into \int_{t_1}^{t_2}g(u)\,du.
    \end{equation*}

    Now from the fact that the Riemann measure is invariant under translation
    follows that also the integral is invariant under translation. Using this
    fact and continuity of the function $f$ we get
    \begin{equation*}
    \begin{split}
        \int_{t_1}^{t_2}g_n(u)\,du &=
        (h_n)^{-1}\left[\int_{t_1}^{t_2}f(u+h_n)\,du-\int_{t_1}^{t_2}f(u)\,du\right]\\
        &=(h_n)^{-1}\left[\int_{t_1+h_n}^{t_2+h_n}f(u)\,du-\int_{t_1}^{t_2}f(u)\,du\right]\\
        &=(h_n)^{-1}\int_{t_2}^{t_2+h_n}f(u)\,du-(h_n)^{-1}\int_{t_1}^{t_1+h_n}f(u)\,du
    \end{split}
    \end{equation*}
    Hence
    \begin{equation*}
         \int_{t_1}^{t_2}g_n(u)\,du\into f(t_2)-f(t_1)
    \end{equation*}
    Therefore we must have
    \begin{equation*}
        \int_{t_1}^{t_2}g(u)\,du= f(t_2)-f(t_1)\fa t_1,t_2\in J.
    \end{equation*}
\ep

\bigskip

\begin{thm}
Let $J$ be an open interval and $\ Y$ a Banach space.
Assume that $g$ from $J$ into the Banach space $Y$
is locally summable on $J$ with respect to Riemann measure.
If the function $g$ is right-side continuous on the interval $J$
then the function
\begin{equation*}
 f(t)=f(t_0)+\int_{t_0}^{t}g(u)\,du\fa t\in J
\end{equation*}
represents a right-sided antiderivative of $g$ over the interval $J.$
\end{thm}

\bigskip
\bp
    The proof is straightforward and we leave it to the reader.
\ep

\bigskip

\section{Lebesgue points of a function\\ summable with respect to
Riemann measure}

\bigskip


A summable function $f$ with respect to the Riemann measure space
is said to have a Lebesgue point at $t=s$ if
\begin{equation*}
    \lim_{h\into 0}\frac{1}{h}\int_s^{s+h}|f(t)-f(s)|\,dt=0.
\end{equation*}

At a Lebesgue point the function $F(t)=\int_{t_0}^t f(u)\,du$ is
differentiable and $F'(s)=f(s).$

\begin{thm}[Almost all points of a summable function are Lebesgue points]
\label{Almost all points of a summable function are Lebesgue
points} Assume that $(R,V,v)$ is the Riemann measure space and $Y$
a Banach space. Then almost every point of any summable function
$f\in L(v,Y)$ is a Lebesgue point. Therefore for every such a
function $f$ the function
 $$F(t)=\int_{t_0}^t f(u)\,du\fa t\in R$$
is differentiable almost everywhere and we have that $F'=f$ almost
everywhere on $R.$
\end{thm}

For a proof of this powerful theorem see Dunford and Schwartz
\cite[page 217, Theorem 8]{DS1}.

\bigskip


\begin{defin}[Absolute continuity of a vector measure]
\label{Absolute continuity of a vector measure} Given a Lebesgue
measure space $(X,V,v)$ on a $\sigma$-algebra $V$ and
$\sigma$-additive vector measure $\mu$  from $V$ into a Banach
space $Y.$ We say that the vector measure $\mu$ is absolutely
continuous with respect to the measure $v$ if \bb
    \mu(A)=0\qtext{whenever} v(A)=0.
\ee
\end{defin}

\bigskip

\begin{thm}[Phillips]
\label{Phillips} Assume that $(X,V,v)$ is a Lebesgue measure space
on a $\sigma$-algebra $V$ such that $v(X)<\infty$ and $Y$ is a
reflexive Banach space.

Assume that $\mu$ is a $\sigma$-additive vector measure of finite
variation from the $\sigma$-algebra $V$ into $Y.$

If $\mu$ is absolutely continuous with respect to the Lebesgue
measure $v,$ then there exists a Bochner summable function $g\in
L(v,Y)$ such that \bb
    \mu(A)=\int_Ag\,dv\fa A\in V.
\ee
\end{thm}

\bigskip

For the proof of this remarkable theorem see Diestel and Uhl
\cite[page 76]{diestel-uhl}. This result can be found in the
original paper of Phillips \cite{phillips}.

\bigskip

\begin{thm}[Integral representation of Lipschitzian functions]
\label{Integral representation of Lipschitzian functions} Let $I$
denote any closed bounded interval of the space $R$ of reals and
$H$ a Hilbert space over the field $R.$ Assume that $(I,V,v)$ is
the Riemann measure space as above. Let $f$ denote a function from
the interval $I$ into the Hilbert space $H.$

If the function $f$ is Lipschitzian on the interval $I,$ that is
\begin{equation*}
    |f(t_1)-f(t_2)|\le m\,|t_1-t_2|\fa t_1,t_2\in I,
\end{equation*}
then there exists a Bochner summable function $g\in L(v,H)$ such
that $\abs{g(t)}\le m$ for all $t\in I$ and
\begin{equation}
    f(t)-f(t_0)=\int_{t_0}^tg(u)\,du\fa t,t_0\in I.
\end{equation}
Moreover the function $f$ is differentiable almost everywhere and
$f'(t)=g(t)$ for almost all $t\in I.$
\end{thm}

\bigskip

\bp
    It follows from Riesz representation theorem that every Hilbert space is reflexive.
    Define a set function $\mu$ on an interval $A\subset I$ with end points $a\le b$
    by
    \bb
        \mu(A)=f(b)-f(a)\qtext{for all intervals} A\in V.
    \ee
    The function $\mu$ is additive on the prering $V.$ Thus it forms a vector measure.
    Since $f$ is Lipschitzian with constant $m$ we have
    \bb
        |\mu(A)|=|f(b)-f(a)|\le m\,(b-a)=m\,v(A)\fa A\in V,
    \ee
    thus the vector measure $\mu$ belongs to the space $K(v,H)$ and $\norm{\mu}\le m.$
    Introduce the bilinear transformation by $u(r,h)=rh$ for all $r\in R$ and $h\in H.$
    Define the set function $\mu_c$ by
    \bb
        \mu_c(A)=\int u(c_A,\mu)\fa A\in V_c.
    \ee
    Since the interval $I$ belongs to the prering $V,$ the family
    $V_c$ of summable sets forms a $\sigma$-algebra. It follows
    form linearity of the integral $\int u(f,d\mu)$ with respect to $f$
    that the set function $\mu_c$ is additive and we have
    \bb
        |\mu_c(A)|=|\int u(c_A,d\mu)|\le |u|\,\norm{c_A}\norm{\mu}=mv_c(A)\fa A\in V_c,
    \ee
    The above implies that the vector measure $\mu_c$ is $\sigma$-additive
    and absolutely continuous with respect to the Lebesgue measure $v_c.$
    Thus by theorem of Phillips \ref{Phillips} there exists a summable
    function $g\in L(v_c,H)=L(v,H)$ such that
    \bb
        \mu_c(A)=\int_A g\,dv_c\fa A\in V_c.
    \ee
    In particular
    \bb
        \mu(A)=\mu_c(A)=\int_Ag\,dv_c=\int_Ag\,dv\fa A\in V
    \ee
    that is
    \bb
        f(t_2)-f(t_1)=\int_{\langle t_1,t_2\rangle}g\,dv
        =\int_{t_1}^{t_2}g(s)\,ds\fa t_1\le t_2,\ t_1,t_2\in I.
    \ee
    It follows from Theorem \ref{When $L(v_2,Y)$ extends $L(v_1,Y)$?}
    that the function $c_Ig$ is Bochner integrable with respect to Riemann measure
    space $(R,V,v).$ Thus almost all points of $c_Ig$ are its Lebesgue points.
    If $s$ lies inside the interval $I$ and $s$
    is a Lebesgue point of the function $c_Ig$ then for sufficiently small $\delta$
    \bb
        |\frac{1}{h}(f(s+h)-f(s))|\le
        \left|\frac{1}{h}\int_s^{s+h}|g(t)-g(s)|\,dt\right|\le m\fa |h|<\delta.
    \ee
    Passing to the limit in the above inequality we get $|f'(s)|=|g(s)|\le m.$
    Replacing eventually the values of $g(s)$ by zero, on a set of measure zero,
    we may assume that $|g(t)|\le m$ for all $t\in I.$
\ep
For a direct proof of the above theorem see Bogdanowicz and Kritt \cite{bogdan22}.

\bigskip
%
%
%
%
%
%

\bigskip

\begin{defin}[Space $M_\infty(I,Y)$]
\label{space of bounded summable functions} Let $(R,V,v)$ be the
Riemann measure space and $I$ a closed bounded interval. Assume
that $Y$ is a Banach space. Denote by $M_\infty(I,Y)$ the set of
all functions that are Bochner summable on the interval $I$ and
bounded almost everywhere. Introduce the functional \bb
    \norm{f}=\inf\set{M:\ |f(x)|\le M \qtext{for almost all}x\in I}.
\ee The functional $\norm{\ }$ forms a seminorm on
$M_\infty(I,Y).$ Identify functions that are equal almost
everywhere. After such identification the pair
$(M_\infty(I,Y),\norm{\ })$ becomes a Banach space. The space
obtained in this way is called the {\bf space of essentially
bounded measurable functions} and the norm {\bf essential sup
norm.}
\end{defin}

\bigskip

\begin{pro}[The space $C(I,Y)$ forms a subspace of $M_\infty(I,Y)$]
\label{The space C forms a subspace of M sub-infinity} The map
$f\mapsto [f],$ where $[f]$ denotes the collection of all
functions that are equal almost everywhere to the function $f,$
establishes an isometry and isomorphism between the space $C(I,Y)$
and a subspace of the space $M_\infty(I,Y).$
\end{pro}


\bp
    The proof is straightforward and we leave it to the reader.
\ep



\bigskip

\end{document}